\documentclass[12pt,leqno]
{amsart}

\usepackage[margin=1.2in]{geometry}

\usepackage[colorlinks, linkcolor=blue,citecolor=blue]{hyperref}

\usepackage{amsmath,epsfig,graphicx,color, float}
\usepackage{subcaption}
\usepackage{relsize}
\usepackage{comment}
\usepackage{algorithm}
\usepackage{algorithmicx}
\usepackage{algpseudocode}
\usepackage{natbib}
\usepackage{multirow}

\numberwithin{table}{section}
\numberwithin{figure}{section}
\numberwithin{equation}{section}

\definecolor{darkblue}{rgb}{.2, 0.2,.8}
\definecolor{darkgreen}{rgb}{0,0.5,0.3}
\definecolor{darkred}{rgb}{.8, .1,.1}

\newcommand{\bfY}{\vect{Y}}
\newcommand{\bfy}{\vect{y}}
\newcommand{\bfx}{\vect{x}}
\newcommand{\bfX}{\vect{X}}

\newcommand{\bfTheta}{\mat{\Theta}}
\newcommand{\bfpi}{\vect{\pi}}
\newcommand{\bfT}{\mat{T}}
\newcommand{\bft}{\vect{t}}

\newcommand{\bfc}{\vect{c}}

\newcommand{\bfbeta}{\vect{\beta}}
\newcommand{\bfeta}{\vect{\eta}}
\newcommand{\bftheta}{\vect{\theta}}

\newcommand{\1}{{\mathbf 1}}
\newcommand{\0}{\mat{0}}

\newcommand{\N}{\mathbb{N}}
\newcommand{\R}{\mathbb{R}}
\newcommand{\E}{\mathbb{E}}
\renewcommand{\P }{{\mathbb P}}

\newtheorem{lemma}{Lemma}[section]

\newtheorem{proposition}[lemma]{Proposition}

\newtheorem{example}[lemma]{Example}

\newtheorem{remark}{Remark}[section]

\usepackage{bm}
\newcommand{\vect}[1]{\pmb{#1}}
\newcommand{\mat}[1]{\boldsymbol{\bm #1}}

\allowdisplaybreaks

\newcommand{\indep}{\perp \!\!\! \perp}

\usepackage{tikz}
\usetikzlibrary{arrows,positioning,calc} 
\usetikzlibrary{decorations.pathreplacing}
\tikzset{
    >=stealth',
    punkt/.style={
           rectangle,
           rounded corners,
           draw=black, thick,
           text width=5.5em,
           minimum height=2em,
           text centered},
    punktl/.style={
           rectangle,
           rounded corners,
           draw=black, thick,
           text width=7em,
           minimum height=2em,
           text centered},
    pil/.style={
           ->,
           shorten <=4pt,
       shorten >=4pt
    },
    pildotted/.style={
           ->,
           shorten <=4pt,
           shorten >=4pt,
  dotted,
  },
    punktf/.style={
           rectangle,
           text width=4.0em,
           minimum height=1.5em,
           text centered},
    punktfTop/.style={
           rectangle,
           text width=4.0em,
           minimum height=1.5em,
           text centered,
           append after command={
               [thick,shorten >=0.2bp, shorten <=0.2bp]
               (\tikzlastnode.north west)edge(\tikzlastnode.north east)
}
    },
    punktfBot/.style={
           rectangle,
           text width=4.0em,
           minimum height=1.5em,
           text centered,
           append after command={
               [thick,shorten >=0.2bp, shorten <=0.2bp]
               (\tikzlastnode.south west)edge(\tikzlastnode.south east)
            }
    }
}

\usepackage[font=footnotesize,labelfont=bf]{caption}

\begin{document}
\title[Bivariate phase-type distributions for experience rating]{Bivariate phase-type distributions for experience rating in disability insurance}
\thanks{}

\author[C.~Furrer]{Christian Furrer}
\address{Department of Mathematical Sciences,
University of Copenhagen,
Universitetsparken 5,
DK-2100 Copenhagen,
Denmark}
\email{furrer@math.ku.dk}

\author[J.J.~S\o rensen]{Jacob Juhl S\o rensen}
\address{AkademikerPension,
Smakkedalen 8,
DK-2820 Gentofte,
Denmark}
\email{jjs@akademikerpension.dk}

\author[J.~Yslas]{Jorge Yslas}
\address{Institute for Financial and Actuarial Mathematics,
University of Liverpool,
L69 7ZL,
Liverpool,
United Kingdom}
\email{Jorge.Yslas-Altamirano@liverpool.ac.uk}

\begin{abstract}
In this paper, we consider the problem of experience rating within the classic Markov chain life insurance framework. 
{
We begin by establishing a link between mixed Poisson distributions and the problem of pricing group disability insurance contracts that exhibit heterogeneity. We focus on shrinkage estimation of disability and recovery rates, taking into account sampling effects such as right-censoring. We then investigate some specific multivariate mixed Poisson models with mixing distributions encompassing independent Gamma, hierarchical Gamma, and multivariate phase-type. In particular, we demonstrate how maximum likelihood estimation for these models can be performed using expectation-maximization algorithms, which might be of independent interest.
}
Finally, we showcase the practicality of the proposed shrinkage estimators through a numerical study based on simulated yet realistic insurance data. Our findings highlight that by allowing for dependency between latent group effects, estimates of recovery and disability rates mutually improve, leading to enhanced predictive performance.

\end{abstract}
\maketitle

\section{Introduction}

This paper builds on and extends the work of~\cite{furrer2019experience} on group shrinkage estimation of transition rates in Markov chains by allowing for dependence between latent group effects across transitions. Following closely~\cite{yslas2021fitting}, but focusing on mixed Poisson regressions rather than frailty models, we study multivariate phase-type distributions as a tractable and flexible class of priors for the latent group effects. The connection to Markov chain models is provided by the classic connection between Poisson regression and parametric inference in multi-state models. The phase-type approach is compared with other alternatives, including the hierarchical system of finite Gamma mixtures proposed in~\cite{norberg1989class} as well as some naïve choices, through an in-depth simulation study related to experience rating for disability insurance. {Overall}, the phase-type approach outperforms the alternatives.

In actuarial science, but also beyond, mixed Poisson regressions constitute an important and highly successful model class for frequency risk, not least due to the fact that, contrary to ordinary Poisson regressions, mixed Poisson regressions may capture overdispersion caused by heterogenity across individuals or groups of individuals. For an introduction with an actuarial flair, see Chapter~2 in~\cite{Denuit2007}. The identification of tractable classes of priors that allow for dependence as well as associated maximum likelihood estimation based on expectation-maximization algorithms~\citep{Dempster1977} has received significant attention; in the following, we just mention a few illustrative contributions. To allow for positive dependence, the simplest would be to assume that the marginal Poisson distributions share a common mixing variable; for this case,~\cite{ghitany2012algorithm} follows~\cite{Karlis2001,Karlis2005} and proposes an expectation-maximization algorithm. In~\cite{Chen2024}, expectation-maximization algorithms for the bivariate lognormal, see also~\cite{AitchisonHo1989} and~\cite{JeongTzougasFung2023}, and the Gaussian copula paired with Gamma marginals are provided; these classes of priors have already been applied to motor third-party liability insurance in~\cite{Pechon2018}. 

Our contribution to this area is two-fold. First, we extend the work of~\cite{norberg1989class} on a hierarchical class of Gamma mixtures to the regression setting. In particular, we develop an expectation conditional maximization algorithm. The downside of this class of priors is that it is only able to capture positive correlation. Second, we adapt the work of~\cite{yslas2021fitting} on phase-type frailty models to the mixed Poisson setting and describe the corresponding expectation-maximization algorithm. In particular, we offer some numerical considerations that are important for an efficient implementation of the algorithm. The underlying idea is that multivariate phase-type priors could offer a more flexible and mathematically tractable alternative, capable of capturing a wider range of dependence structures.

In its univariate form, a phase-type distribution is defined as the distribution of the time until absorption of a time-homogeneous, finite state-space Markov jump process with a single absorbing state, while the remaining states are transient. These distributions have found several applications in applied probability due to closed-form formulas for different functionals and, more importantly, their denseness on the set of distributions on the positive real line. In actuarial science, we find applications towards areas such as risk theory~\citep{Bladt2005}, loss severity modeling~\citep{BladtYslas2023}, mortality modeling~\citep{AlbrecherBladtBladtYslas2022}, and multi-state modeling~\citep{AhmadBladtFurrer2023,AhmadBladt2023}. Several multivariate extensions of phase-type distributions exist in the literature. Notable among them are the general MPH* class introduced in~\cite{kulkarni1989new} and its more computationally amenable subclasses: the feed-forward type class, confer, for example, with~\cite{albrecher2020fitting}, and the mPH class recently detailed in~\cite{bladt2023tractable}, which also includes applications towards severity modeling. In this work, we opt as in~\cite{yslas2021fitting} for the feed-forward type class due to its mathematical convenience and its denseness on the set of the distributions on the positive orthant.

Both classes of priors are compared to some naïve alternatives through an in-depth simulation study on group disability insurance. {For disability insurance, there is often a natural group structure, especially if the insurance is provided through a corporate (pension) scheme, in which case each company constitutes a more or less closed group. Due to differences in characteristics across companies, including differences in internal practices, there is often substantial heterogeneity between groups -- also compared to the heterogeneity within groups -- which in turn gives rise to non-negligible group effects. The} numerical study serves a dual purpose: It demonstrates that our new specifications can be implemented smoothly for near real-word scenarios and it provides insight into the advantages and disadvantages of each specification. The phase-type specification {essentially} exhibits as good or strictly better predictive performance across all relevant cases.

For the application to group disability insurance, we draw on the work of~\cite{furrer2019experience}, which employed an empirical Bayes approach to group experience rating in Markov chain models with the latent effects introduced on the level of transition rates. {Linear Bayes methodology, rather than empirical Bayes, constitutes an alternative which in actuarial science is often referred to as credibility theory; for an overview, see~\cite{buhlmanngisler2005}.} Although~\cite{furrer2019experience} provides details on the sufficiency of summary statistics such as occurrences and exposures, hereby replicating the by now well-known link between likelihood based inference for multi-state and Poisson models, the proposed shrinkage estimators are derived under a strict assumption of independence between the latent group effects. However, dependence between latent effects is to be expected, especially if only few covariates are included. Thus, allowing for dependence could enhance predictive accuracy. Concordantly, in~\cite{furrer2019experience} such an extension is highlighted as a promising avenue for future research. In light of this, the present work may be viewed as picking up exactly where~\cite{furrer2019experience} left, showcasing the potential of a more sophisticated modeling approach.

The remainder of the paper is organized as follows. 
{
In Section~\ref{sec:experience}, we describe the problem of experience rating for disability insurance based on Markov chain modeling and recall the connection to mixed Poisson regression. In Section~\ref{sec:mpr}, we provide an overview of mixed Poisson regression and thoughtfully analyze the model's properties and its estimation when the mixing distribution is one of three distinct mixing distributions: independent Gamma, hierarchical Gamma, and multivariate phase-type. 
}
Section~\ref{sec:numerical} is dedicated to a detailed numerical study based on simulated data showcasing the practical applicability of the presented models as well as the overall preferability of the phase-type specification. Finally, Section~\ref{sec:conclusion} concludes. 

\section{Experience rating for disability insurance}\label{sec:experience}

{The area of application of this paper is experience rating in disability insurance; that is, we seek to price group disability insurance contracts subject to heterogeneity}. Subsection~\ref{sec:disability_coverage} contains a description of the insurance coverage of interest. In Subsection~\ref{sec:disability_estimation}, we consider shrinkage estimation of disability and recovery rates subject to, among other things, right-censoring and following closely~\cite{furrer2019experience}. In particular, we provide a link to {mixed Poisson models, which are the focal point of Section~\ref{sec:mpr}.}

\subsection{Disability coverage}\label{sec:disability_coverage}

We specify the insurance coverage by stipulating the benefits $B = (B(t))_{t\geq0}$ paid by the insurer to the insured. We let
\begin{align}\label{eq:disability_payments}
B(\mathrm{d}t)
=
1\{t - U(t) \leq \tau\} \Big(
1\{\varepsilon \leq U(t)\} 1\{Z_t = i\} b_i(t) \, \mathrm{d}t
+ b_{ai}(t) N_{ai}(\mathrm{d}t)
\Big)
\end{align}
with $B(0) = 0$. Here $Z=(Z_t)_{t\geq0}$ denotes the state of the insured, taking values in the set $\mathcal{S} = \{a, i, d\}$ with state $a$ denoting `active', state $i$ denoting `invalid' (disabled), and state $d$ denoting `dead'. Furthermore, $N_{jk}=(N_{jk}(t))_{t\geq0}$ is the counting process associated with the transitions from state $j$ to state $k$, so that $N_{jk}(t)$ denotes the number of transitions from state $j$ to state $k$ in the time interval $[0,t]$, and $U=(U(t))_{t\geq0}$ is the duration since the last transition given by
\begin{align*}
U(t) = t - \sup\{ s \in [0,t] : Z_s \neq Z_t\}.
\end{align*}
Finally, $\varepsilon \in [0,\infty)$ is a waiting period, $\tau \in (0,\infty]$ is the coverage period, $t \mapsto b_{ai}(t)$ are deterministic benefits upon disability, and $t \mapsto b_i(t)$ is a deterministic benefit rate during disability. {The interpretation of~\eqref{eq:disability_payments} is as follows. If the insured becomes disabled in the coverage period, say at time $t$ with $t \leq \tau$, then they receive at exactly time $t$ the payment $b_{ai}(t)$. Additionally, they receive a payment rate of $s \mapsto b_i(s)$ from time $t+\varepsilon$ until recovery or death, whichever occurs first; should it happen that they recover or die before time $t + \varepsilon$, they receive however only the payment $b_{ai}(t)$.}
\begin{example}\label{eq:annuity}
The classic (continuous) disability annuity of rate $b>0$ is recovered by taking $\varepsilon = 0$, $\tau = \infty$, $b_i(t) = b$, and $b_{ai}(t) = 0$, when $B(\mathrm{d}t) = 1\{Z_t = i\} b \, \mathrm{d}t$. In comparison, the general expression~\eqref{eq:disability_payments} also includes a payment upon disability, a coverage period, and a waiting period.
\end{example}
In the following, we suppose that the payment functions $b_{ai}$ and $b_i$ are suitably regular{, for instance piecewise continuous,} and that there exists a maximal contract time $T \in (0,\infty)$ such that $b_{ai}(t) = b_i(t) = 0$ for $t > T$. This is always the case in practice (with $T$ related to, for instance, the legal retirement age). Furthermore, we suppose that the initial state of $Z$ is $a$, and we naturally take $d$ to be absorbing. 

We further assume that $Z$ is Markovian and admits transition rates $t \mapsto \mu_{jk}(t)$, $j, k \in \mathcal{S}$, $j \neq k$, such that
\begin{align*}
\mu_{jk}(t) = \lim_{h\downarrow 0} \frac{\mathbb{P}(Z_{t+h} = k \, | \, Z_t = j)}{h}
\end{align*}
almost everywhere. Since state $d$ is absorbing, we have $\mu_{da}(t) = \mu_{di}(t) = 0$.

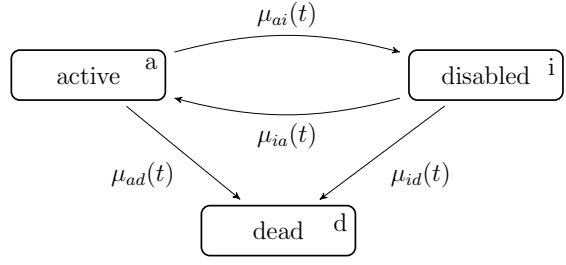
\begin{figure}[!htbp]
	\centering
	\scalebox{0.8}{
	\begin{tikzpicture}[node distance=2em and 0em]
		\node[punkt] (1) {disabled};
		\node[anchor=north east, at=(1.north east)]{i};
		\node[punkt, left = 40mm of 1] (0) {active};
		\node[anchor=north east, at=(0.north east)]{a};
		\node[draw = none, fill = none, left = 20 mm of 1] (test) {};
		\node[punkt, below = 20mm of test] (2) {dead};
		\node[anchor=north east, at=(2.north east)]{d};
	\path
		(0)	edge [pil, bend left = 15]		node [above]				{$\mu_{ai}(t)$}				(1)
		(1)	edge [pil]					node [below right]		{$\mu_{id}(t)$}				(2)
		(0)	edge [pil]					node [below left]			{$\mu_{ad}(t)$}				(2)
		(1)	edge [pil, bend left = 15]		node [below]				{$\mu_{ia}(t)$}				(0)
	;
	\end{tikzpicture}}
	\caption{The Markov chain $Z$ with values in the set $\mathcal{S}=\{a,i,d\}$ and transition rates $t \mapsto \mu_{jk}(t)$.}
	\label{fig:multi_state}
\end{figure}

In other words, we have
\begin{align*}
\mathbb{E}\big(N_{jk}(\mathrm{d}t) \, | \, (Z_s)_{0\leq s<t} \big) = \mathbb{E}\big(N_{jk}(\mathrm{d}t) \, | \, Z_{t-} \big) = \mu_{Z_{t-}k}(t) \, \mathrm{d}t.
\end{align*}
The model is illustrated in Figure~\ref{fig:multi_state}.

For practitioners, an extension to semi-Markov models (where the transition rates may depend on the duration since the last transition) would seem of great interest. However, we should like to stress that the main insights of this paper also carry over to semi-Markov models; actually, this is a rather straightforward extension. To keep the presentation concise and focused on the novelties, we therefore restrict the analysis to Markov chains.

Pricing and reserving based on the equivalence principle necessitates the calculation of the prospective reserves $t \mapsto V_a(t)$ and $(t,s) \mapsto V_i(t,s)$  given by
\begin{align*}
V_a(t) &= \mathbb{E}\bigg(\int_t^T e^{-\int_t^s r(u) \, \mathrm{d}u} B(\mathrm{d}s) \, \bigg| \, Z_t = a\bigg), \\
V_i(t,s) &= \mathbb{E}\bigg(\int_t^T e^{-\int_t^v r(u) \, \mathrm{d}u} B(\mathrm{d}v) \, \bigg| \, Z_t = i, U_t = t - s\bigg),
\end{align*}
where $t \mapsto r(t)$ is a suitably regular{, say piecewise continuous,} interest rate of choice. We are, of course, implicitly requiring that $s \leq t$. The following result is a consequence of, for instance, Corollary~7.2 in~\cite{adechristiansen2017thiele}.
\begin{proposition}\label{prop:thiele}
It holds almost everywhere that
\begin{align}
\begin{split}\label{eq:reserve_comp}
\frac{\mathrm{d}}{\mathrm{d}t}V_a(t) &= \big( r(t) + \mu_{ad}(t)\big) V_a(t) - 1\{t \leq \tau\} \big( b_{ai}(t) + V_i(t,t) - V_a(t) \big) \mu_{ai}(t), \\
\frac{\mathrm{d}}{\mathrm{d}t}V_i(t,s) &= \big( r(t) + \mu_{id}(t)\big) V_i(t,s) \\ &\quad- 1\{s \leq \tau\}\Big(1\{\varepsilon \leq t-s\} b_i(t) + \big( V_a(t) - V_i(t,s)\big) \mu_{ia}(t) \Big).
\end{split}
\end{align}
\end{proposition}
Section~8 in~\cite{adechristiansen2017thiele} provides a numerical scheme for the prospective reserves based on~\eqref{eq:reserve_comp} in conjunction with the obvious boundary conditions $V_a(T) = 0$ and $V_i(T,s) = 0${, the system given by~\eqref{eq:reserve_comp} being a special case of Thiele's differential (or integral) equations for semi-Markov models. The interpretation of~\eqref{eq:reserve_comp} is similar to the classic case with deterministic sojourn payment rates and transition payments, confer with~\cite{Hoem1969}.} 
\begin{remark}
If there is no waiting period, corresponding to $\varepsilon = 0$, then $V_i(t,s) = 1\{s \leq \tau \}W_i(t)$, where
\begin{align*}
\frac{\mathrm{d}}{\mathrm{d}t}W_i(t) &= \big( r(t) + \mu_{id}(t)\big) W_i(t) - b_i(t) -  \big(V_a(t) - W_i(t)\big) \mu_{ia}(t), \quad W_i(T) = 0.
\end{align*}
This leads to a particularly straightforward numerical scheme.
\end{remark}

\subsection{Shrinkage estimation}\label{sec:disability_estimation}

In general, the transition rates $t \mapsto \mu_{jk}(t)$ are unknown and must be estimated from data. This is the trademark of the multi-state approach to modeling disability claims, as opposed to a more direct modeling of the disability's frequency and severity. In disability insurance, many claims remain open at the time of estimation. Expressed differently, the observations are subject to right-censoring. A major advantage of the multi-state approach is the transparent and rigorous manner in which it deals with, for instance, right-censoring. In the context of experience rating, on the other hand, a major drawback of the multi-state approach is its discord with classic credibility theory, so-called linear Bayes theory, confer with the discussion in Section~1 of~\cite{furrer2019experience}.

In the following, we suppose for simplicity that the mortality rates are known. We are thus interested in estimating the disability rate $t \mapsto  \mu_{ai}(t)$ and the recovery rate $t \mapsto \mu_{ia}(t)$ on $[0,T]$. To allow for random entries and exits in the insurance portfolio, we introduce a filtering process $C=(C(t))_{t\geq0}$, so that $C(t) = 1$ if the insured is in the portfolio at time $t$ and $C(t) = 0$ if the insured is not in the portfolio at time $t$. The case of only right-censoring corresponds to $C(t) = 1\{t \leq R\}$ for some right-censoring time $R$. Under suitable assumptions, including independent filtering, the relevant partial likelihood process $\mathcal{L}$ reads $\mathcal{L} = \mathcal{L}_{ai} \mathcal{L}_{ia}$ with
\begin{align*}
\mathcal{L}_{ai} &= \exp\bigg\{\int_0^T C(t) \log (\mu_{ai}(t)) \, N_{ai}(\mathrm{d}t) - \int_0^T C(t) 1\{Z_{t-} = a\} \, \mathrm{d}t\bigg\}, \\
\mathcal{L}_{ia} &= \exp\bigg\{\int_0^T C(t) \log (\mu_{ia}(t)) \, N_{ia}(\mathrm{d}t) - \int_0^T C(t) 1\{Z_{t-} = i\} \, \mathrm{d}t\bigg\},
\end{align*}
see Section~III.4 in~\cite{andersen1993counting}. In actuarial practice, estimation of the disability and recovery rates is often based on Poisson distributions. The validity of such an approach rests on the likelihood principle and a simple approximation, namely that the disability rate $t \mapsto \mu_{ai}(t)$ and recovery rate $t \mapsto \mu_{ia}(t)$ are continuous from the left and piecewise constant on some grid $0 = t_0 < t_1 < \cdots < t_K = T$, $K \in \mathbb{N}$, of $[0,T]$. In this case, it holds that
\begin{align} \label{eq:L_ai}
\mathcal{L}_{ai} &= \prod_{k=1}^K {(\mu_{ai}(t_k))}^{O_{ai}^k} \exp(-E_a^k \mu_{ai}(t_k)), \\ \label{eq:L_ia}
\mathcal{L}_{ia} &= \prod_{k=1}^K {(\mu_{ia}(t_k))}^{O_{ia}^k} \exp(-E_i^k \mu_{ia}(t_k)), 
\end{align}
where $(E_a^k)_{k=1}^K$, $(E_i^k)_{k=1}^K$, $(O_{ai}^k)_{k=1}^K$, and $(O_{ia}^k)_{k=1}^K$ are the relevant exposures and occurrences given by
\begin{align*}
&E_a^k = \int_{t_{k-1}}^{t_k} C(t) 1\{Z_{t-} = a\} \, \mathrm{d}t, &&E_i^k = \int_{t_{k-1}}^{t_k} C(t) 1\{Z_{t-} = i\} \, \mathrm{d}t, \\
&O_{ai}^k = \int_{t_{k-1}}^{t_k} C(t) N_{ai}(\mathrm{d}t), &&O_{ia}^k = \int_{t_{k-1}}^{t_k} C(t) N_{ia}(\mathrm{d}t).
\end{align*}
It follows that the likelihoods~\eqref{eq:L_ai} and~\eqref{eq:L_ia} are proportional to the density functions that appear by considering the exposures fixed and known and the occurrences independent with the following marginal distributions:
\begin{align*}
O_{ai}^k \sim \text{Poisson}\big(E_a^k \mu_{ai}(t_k)\big), \quad \quad O_{ia}^k \sim \text{Poisson}\big(E_i^k \mu_{ia}(t_k)\big).
\end{align*}
This establishes the desired link to Poisson distributions. In general, links of this type are well-known in the literature, confer with the discussion in Subsection~2.2 of~\cite{furrer2019experience}.

Experience rating for group disability insurance {is closely related to shrinkage estimation of group effects. Experience rating concerns the improvement of a portfolio level assessment through the use of information about a specific (sub-)group's performance. Shrinkage estimators improve, either explicitly or implicitly, a naïve estimate by bringing it closer to a value supplied by additional information. If we think of the raw group effect as the naïve estimate and the portfolio level as the value supplied by additional information, we see that shrinkage estimation and experience rating seek the same goal, be it from different directions.}

{In this context},~\cite{furrer2019experience} suggests estimation of the transition rates using empirical Bayes methods. Disregarding mortality, the `simple shrinkage' estimation of~\cite{furrer2019experience} is based on the following structural assumptions:
\begin{align} \nonumber
&\mathbb{E}\big(N_{ai}(\mathrm{d}t) \, \big| \, Z_{t-} = a, (\Theta_{ai},\Theta_{ia}) = (\theta_{ai}, \theta_{ia}) \big)
=
\theta_{ai} \mu_{ai}(t) \, \mathrm{d}t, \\ \nonumber
&\mathbb{E}\big(N_{ia}(\mathrm{d}t) \, \big| \, Z_{t-} = i, (\Theta_{ai},\Theta_{ia}) = (\theta_{ai}, \theta_{ia}) \big)
=
\theta_{ia} \mu_{ia}(t) \, \mathrm{d}t, \\ \label{eq:gamma_indep_furrer2019}
&\Theta_{ai} \indep \Theta_{ia}, \quad \Theta_{ai} \sim \text{Gamma}(\psi_{ai}^{-1},\psi_{ai}^{-1}), \quad \Theta_{ia} \sim \text{Gamma}(\psi_{ia}^{-1},\psi_{ia}^{-1}),
\end{align}
with $(\Theta_{ai}, \Theta_{ia})$ being latent group effects. The transition rates are identical between groups, while the latent effects differ; we have not highlighted this in the notation. Regarding estimation, the introduction of latent group effects just requires replacing the Poisson distributions with mixed Poisson distributions. This is under the unverifiable assumption that the censoring mechanism is conditionally non-informative of the latent effects, see Subsection~IX.3 in~\cite{andersen1993counting}. In other words, the likelihoods of interest correspond to the following mixed Poisson models:
\begin{align*}
O_{ai}^k \, | \, (\Theta_{ai},\Theta_{ia}) { = (\theta_{ai},\theta_{ia})} &\sim \text{Poisson} \big({\theta_{ai}} E_a^k \mu_{ai}(t_k)\big), \\
O_{ia}^k \, | \, (\Theta_{ai},\Theta_{ia}) { = (\theta_{ai},\theta_{ia})} &\sim \text{Poisson} \big({\theta_{ia}} E_i^k \mu_{ia}(t_k)\big).
\end{align*}
This provides a direct link to Section~\ref{sec:mpr}.

The `simple shrinkage' estimation of~\cite{furrer2019experience} operates under assumptions~\eqref{eq:gamma_indep_furrer2019}, that is, independent Gamma mixing, and produces empirical Bayes estimates of the disability and recovery rates.  The disability coverage may then be priced by plugging these estimates into~\eqref{eq:reserve_comp}. Unfortunately, the prediction errors resulting from such a `plug-in' approach are generally non-traceable. It is, therefore, important to assess the final predictive performance carefully using relevant measures of performance. 

The independence assumption is disadvantageous, as in practice you typically observe substantial dependence between the group effects, especially if the group effects indirectly capture differences in the types of disabilities experienced by individuals across groups. Stress related claims are, for example typically more frequent but also shorter in duration than claims stemming from diseases of or injuries to the musculoskeletal system. This entails positive correlation between the disability and recovery rate. In light hereof, {the next section contains a theoretical investigation of mixed Poisson models with other types of mixing than that of~\eqref{eq:gamma_indep_furrer2019}, namely hierarchical and phase-type mixings. This theoretical investigation is finally supplemented by a numerical study devoted to experience rating for disability insurance. By allowing for dependence between latent group effects, estimation of the recovery rate may borrow strength from the estimation of the disability rate -- and vice versa. Everything else being equal, this should lead to improvements in predictive performance.}

\section{Mixed Poisson regression}\label{sec:mpr}

Mixed Poisson regression is a versatile tool for modeling count data, especially when observations exhibit overdispersion or correlation. By introducing random effects into the Poisson regression framework, the mixed Poisson model can capture unobserved heterogeneity and allow for more complex dependence structures between observations. This flexibility makes it suitable for a wide range of applications, including insurance counts modeling. In this section, we delve into the specificities of this model in the multivariate case, with a particular emphasis on three tractable mixing distributions: independent Gamma (Subsection~\ref{subsec:ind}), hierarchical Gamma (Subsection~\ref{subsec:hier}), and multivariate phase-type (Subsection~\ref{subsec:ph}).

\subsection{Setup}

Consider a random vector $\bfY = (Y_1, \dots, Y_d)$ with conditionally independent Poisson-distributed margins given a nonnegative random variable $\Theta$, known as the mixing component. More specifically, we assume that the laws of the conditional margins are 
\begin{align*}
	Y_i \mid  \Theta = \theta \sim \text{Poisson}(\theta E_i e^{\bfX_i \bfbeta}) ,
\end{align*}
where $\bfbeta$ is an $h$-dimensional column vector of regression coefficients, $\bfX_i$ are $h$-dimensional row vectors of covariates, and $E_i$ are exposure factors, $i = 1, \dots, d$. In this way, the conditional joint density function $f_{\bfY | \Theta}$ of $\bfY$ given $\Theta$ reads 
\begin{align*}
	f_{\bfY | \Theta}(\bfy | \theta) = \prod_{i = 1 }^{d} \frac{(\theta E_i e^{\bfX_i \bfbeta})^{y_i}}{y_i!}\exp(-\theta E_ie^{\bfX_i \bfbeta}), \quad \bfy \in \N_0^d. 
\end{align*}
Thus, the general form of the (unconditional) joint density $f_{\bfY}$ of $\bfY$ is 
\begin{align} \label{eq:mmpr}
	f_{\bfY}(\bfy) = \int_0^\infty \prod_{i = 1 }^{d} \frac{(\theta E_ie^{\bfX_i \bfbeta})^{y_i}}{y_i!}\exp(-\theta E_i e^{\bfX_i \bfbeta}) f_{\Theta}(\theta)  \, \mathrm{d}\theta,
\end{align}
where $f_{\Theta}$ denotes the density function of $\Theta$. An alternative expression for the joint density of $\bfY$ can be given in terms of the Laplace transform $\mathcal{L}_{\Theta}(u) = \mathbb{E}[\exp(-u\Theta)]$, $u \geq 0$, of $\Theta$ as follows
\begin{align*}
	f_{\bfY}(\bfy) =  \left(\prod_{i = 1 }^{d} \frac{( E_i e^{\bfX_i \bfbeta})^{y_i}}{y_i!}\right)(-1)^{\sum_{i = 1}^d y_i}\mathcal{L}_{\Theta}^{(\sum_{i = 1}^d y_i)}\left( \sum_{i = 1}^d E_ie^{\bfX_i \bfbeta}\right)\!,
\end{align*}
where $\mathcal{L}_{\Theta}^{(n)}$ denotes the derivative of order $n \in \N_0$ of $\mathcal{L}_{\Theta}$, provided they exist. 

\begin{example}[Gamma mixing] \label{ex:gammix} \rm
	Assume that $\Theta \sim \text{Gamma}(\alpha, \kappa)$, $\alpha, \kappa >0$, so that
	\begin{align*}
		 \mathcal{L}_{\Theta}(u) = \left( 1 + u / \kappa \right)^{-\alpha}, \quad u \geq 0 .
	\end{align*}
	Then, the joint density function of $\bfY$ is given by
\begin{align*}
	f_{\bfY}(\bfy)  = \left(\prod_{i = 1 }^{d} \frac{( E_ie^{\bfX_i \bfbeta})^{y_i}}{y_i!}\right)\frac{\Gamma(\sum_{i = 1}^d y_i + \alpha)}{ \Gamma(\alpha)}  \frac{\kappa^\alpha}{\left(\kappa + \sum_{i = 1}^d E_ie^{\bfX_i \bfbeta}\right)^{\sum_{i = 1}^d y_i + \alpha }},
\end{align*} 
which corresponds to the joint density of a multivariate negative binomial distribution. 
To avoid overparametrization, it is customary to assume that $\E(\Theta) = 1$. In the Gamma mixing case, this is accomplished by setting $\alpha = \kappa = \psi^{-1}$, $\psi > 0$, so that $\mathbb{E}(\Theta) = 1	$ and $\text{Var}(\Theta) = \psi$.
\hfill $\circ$
\end{example}

We are now interested in extending the above model to a setting where we have different characteristics (or risks) for various groups of individuals. 
Thus, consider a group of $d$ individuals, each possessing $q$ characteristics of interest. We represent by the random vector $\bfY_j = (Y_{1j}, \dots, Y_{dj})$ the $j$-th characteristic of the group, $j = 1, \dots, q$. Given a random vector $\bfTheta = (\Theta_1,\dots, \Theta_q)$, we assume that the random variables $Y_{ij}$ are mutually independent and Poisson distributed. Concretely, we assume that
\begin{align*}
	Y_{ij} \mid  \bfTheta = \bftheta \sim \text{Poisson}(\theta_j E_{ij}e^{\bfX_{ij} \bfbeta_{j}}),
\end{align*}
where $\bfbeta_{j}$ are $h_j$-dimensional column vectors of regression coefficients, $\bfX_{ij}$ are $h_j$-dimensional row vectors of covariates, and $E_{ij}$ are exposure factors.  Then we have that the (unconditional) joint density of $\bfY = (\bfY_1, \dots, \bfY_q)$ is given by
\begin{align}\label{eq:mmpr_gen}
	f_{\bfY}(\bfy ) = \int_{\R^q_+} \prod_{j = 1 }^{q} \prod_{i = 1 }^{d} \frac{(\theta_j E_{ij}e^{\bfX_{ij} \bfbeta_{j}})^{y_{ij}}}{y_{ij}!}\exp(-\theta_j E_{ij} e^{\bfX_{ij} \bfbeta_{j}}) f_{\bfTheta}(\bftheta)  \, \mathrm{d}\bftheta,
\end{align}
where $f_{\bfTheta}$ denotes the joint density of $\bfTheta$. Note that we have considered the case where the regression coefficients differ by characteristic but not by individual. Nevertheless, individual-specific variation of coefficients is a straightforward extension. 

In the following subsections, we provide an in-depth analysis of three choices of prior distributions of $\bfTheta$ that lead to tractable models that allow for efficient estimation. In the last case, we focus mainly on the case $q = 2$ for clarity of presentation -- but also in accordance with the later application and numerical study. 

\textbf{Notation}: To make the mathematical expressions more transparent and concise, we adopt the following notation throughout the remainder of the paper: $e_{ij} := E_{ij} e^{\bfX_{ij} \bfbeta_{j}}$, $e_{\bullet j} := \sum_{i = 1}^{d} e_{ij}$, and $y_{\bullet j} := \sum_{i = 1}^{d} y_{ij}$. 

\begin{remark}
In all cases, the covariates only enter the model through the conditional distributions of the responses given $\bf\Theta$; the prior distribution of $\bf\Theta$ is not allowed to depend on the covariates. This assumption, although standard for mixture regressions, can be hard to verify and might substantially limit which covariates can be used. Confer with~\cite{BoucherDenuit2006} for analyses in the case where the random effects and covariates are correlated.
\end{remark}

\subsection{Independent Gamma mixing} \label{subsec:ind}

We start by considering a mixing component  $\bfTheta = (\Theta_1, \dots, \Theta_q)$ with $\Theta_j \sim \text{Gamma}(\gamma^{-1}_j, \gamma^{-1}_j)$, $\gamma_j > 0$, $j =1, \dots, q$, independent. Under this assumption, we have that~\eqref{eq:mmpr_gen} takes the explicit form
\begin{align*}
	f_{\bfY}(\bfy)  = \prod_{j = 1 }^{q}  \left(\prod_{i = 1 }^{d} \frac{( e_{ij})^{y_{ij}}}{y_{ij}!}\right) \frac{\Gamma(y_{\bullet j} + \gamma^{-1}_j)\gamma^{-\gamma^{-1}_j}_j}{\Gamma(\gamma^{-1}_j) \left(\gamma^{-1}_j + e_{\bullet j}\right)^{y_{\bullet j} + \gamma^{-1}_j }}   ,
\end{align*} 
meaning that the random vectors $\bfY_j$, $j = 1,\dots, q$, are independent multivariate negative binomial distributed. Furthermore, under this specification, we can employ that the Gamma distribution is a conjugate prior for the Poisson distribution, which facilitates the computation of several quantities of interest, such as Bayes estimators. More specifically, given observed data $\bfy =(\bfy_1, \dots,  \bfy_q) $, the posterior distribution for $\Theta_j$ is again Gamma with the explicit representation 
\begin{align*}
\Theta_j \mid \bfY = \bfy \sim 	\text{Gamma}\!\left(y_{\bullet j} +  \gamma^{-1}_j, e_{\bullet j} +\gamma^{-1}_j\right)\!, \quad  j = 1, \dots, q.
\end{align*}
Regarding estimation, it is clear that the fitting of this model to data reduces to estimating $q$ different multivariate mixed Poisson model. We also see that if there are instead common parameters across characteristics, the only added difficulty is that estimation must be performed simultaneously. Typically, this is done by using an expectation-maximization (EM) algorithm. The need for an EM algorithm arises from the nature of the mixing component, which is unobserved. This leads to an incomplete data setup, a scenario in which the EM algorithm has shown to be particularly effective; see \cite{mclachlan2007algorithm} for a comprehensive account of the EM algorithm. In particular, for the Gamma mixing case, an explicit EM algorithm is detailed in for example~\cite{ghitany2012algorithm}. 

\subsection{Hierarchical Gamma mixing} \label{subsec:hier}

The previous assumption of independence may not be realistic in certain scenarios. Hence, we now relax the independence assumption by considering the class of hierarchical models introduced in~\cite{norberg1989class}. More specifically, we consider  $\bfTheta = (\Theta_1, \dots, \Theta_q)$ such that, conditionally on $\Theta_0 \sim \text{Gamma}(\eta, \nu)$, $\eta, \nu >0$, $\Theta_j$  are independent with laws 
\begin{align*}
	\Theta_j \mid \Theta_0 {=\theta_0} \sim \text{Gamma}({\theta_0}, \delta) , \quad j = 1, \dots,q , \quad\delta >0 .
\end{align*}
To avoid overparametrization, we assume that $\E(\Theta_j) = 1$ for all $j = 1, \dots, q$, which is equivalent to taking $\delta = \eta \nu ^{-1}$. Under this specification, we have that 
\begin{align*}
	\text{Var}(\Theta_j ) = \frac{\nu + 1}{\eta} , \quad \text{Cov}(\Theta_j , \Theta_k) = \frac{1}{\eta} , \quad \text{Corr}(\Theta_j , \Theta_k) = \frac{1}{1 + \nu} ,
\end{align*}
for all $j,k = 1, \dots, q$, $j \neq k$. In particular, from above, we note that the hierarchical framework allows only for positively correlated mixing variables. {Furthermore, the variances of the marginal distributions cannot differ. Therefore, even in the limiting case $\eta \to \infty, \nu \to \infty, \eta / \nu \to \delta$, we do not recover the previous independent Gamma mixing.}
 
Unfortunately, no fully explicit expression for the joint density~\eqref{eq:mmpr_gen} can be obtained under this specification of $\bfTheta$. However, it can be computed numerically using the following representation obtained easily by applying the tower property:
\begin{align*}
	f_{\bfY}(\bfy)  =\left(\prod_{j = 1}^q \prod_{i = 1}^d\frac{ (e_{ij})^{y_{ij}}}{  y_{ij}!} \right)\E\!\left( \frac{\delta^{q\Theta_0} \prod_{j = 1}^q \Gamma(\Theta_0 + y_{\bullet j})}{\Gamma(\Theta_0)^q \prod_{j = 1}^q \left(\delta + e_{\bullet j}\right)^{\Theta_0 + y_{\bullet j}}}  \right)\!.
\end{align*}

Interestingly, this model has an interpretation in terms of conjugate priors. More specifically, given observed data $\bfy$, we have that
 \begin{align}\label{eq:hiermar}
 	\Theta_j \mid \Theta_0 {=\theta_0} , \bfY = \bfy \sim \text{Gamma}\Big(  y_{\bullet j} +{\theta_0}, e_{\bullet j} + \frac{\eta}{\nu} \Big) .
 \end{align}
Moreover, it is easy to see that
\begin{align}\label{eq:hiermix}
 	f_{\Theta_0 |\bfY}(\theta |\bfy) 
 	&\propto  \theta^{\nu - 1}e^{-\eta \theta} \prod_{j = 1}^q \frac{\Gamma(y_{\bullet j} + \theta)}{\Gamma(\theta)}\Big( 1+ (\nu / \eta)e_{\bullet j}\Big)^{-\theta} ,
 \end{align} 
which can be identified as the density function of a mixture of Gamma distributions. To obtain the exact specification of this mixture, define $\vect{o} = (o_0 , \dots, o_r ) \in \N_0^{r+1}$ with  
$o_m = \#\{j : y_{\bullet j} >m \}$ 
and 
$r = \max \{\max_j y_{\bullet j} -1, 1 \}$, let $\vect{0}$ be the vector of zeroes, and let $\bfc_{s}$ be the vector with the $s$'th entry equal to one and the remaining entries equal to zero. Then, the mixture distributions are
\begin{align*}
 	\text{Gamma}\Big(m + \eta, \sum_{j = 1}^{q}\log\big(1 + (\nu / \eta) e_{\bullet j}\big) + \nu \Big),
\end{align*}
with corresponding (unscaled) mixture weights $w_m$ given by 
  \begin{align*}
 	w_m = a_m(\vect{o}) \cdot b_m , \quad m = o_0, \dots, \sum_{j = 1}^q y_{\bullet j},
 \end{align*}
 where the coefficients can be computed via recurrence relations. First, for $(b_m)_m$ one can use
\begin{align*}
b_m &=b_{m-1} \frac{m+\eta-1}{\sum_{j = 1}^{q}\log\big(1 + (\nu/ \eta) e_{\bullet j}\big) + \nu},
\end{align*}
initialized by $b_0 = 1$. Second, for $(a_m(\mathbf{o}^{\star}))_{m,\mathbf{o}^{\star}}$ one may use
\begin{align*}
	a_m\left(\mathbf{o}^{\star}+\bfc_{s}\right)=s \cdot a_m\left(\mathbf{o}^{\star}\right)+a_{m-1}\left(\mathbf{o}^{\star}\right)\!,
\end{align*}
initialized by $a_m(\vect{0}) = 1\{m = 0\}$ and taking $a_m(\cdot) = 0$ for $m \notin \{o_0, \ldots, \sum_{j = 1}^q y_{\bullet j}\}$. For an efficient implementation, one should take a shortest path (across $m$) starting from $\vect{0}$ and finishing at $\textbf{o}$.

The scaled mixing weights $\bar{w}_m$ satisfying $\sum_m \bar{w}_m = 1$ can then be computed straightforwardly from the unscaled mixing weigths $w_m$. This interpretation in terms of conjugate priors is particularly useful for the derivation of the estimation procedure presented next. Practical experiences indicate, unfortunately, that this recurrence scheme might {lead to arithmetic overflow}  for large datasets.

\subsubsection{Estimation}
We now show how maximum likelihood estimation of this model can be carried out via an ECM algorithm. To the best of our knowledge, the derivation of this explicit estimation method has not been presented before. Recall that the ECM algorithm \citep[confer with][]{meng1993maximum} is an extension of the EM algorithm, which replaces the maximization (M) step with several computationally simpler conditional maximization (CM) steps.

Consider $\tilde{\bfy} = \{\bfy^{(1)}, \dots,\bfy^{(N)}\}$: an iid sample of size $N$ coming from the model specified by \eqref{eq:mmpr_gen} with hierarchical Gamma mixing and corresponding paired covariate information $\tilde{\bfx} = \{\bfx^{(1)}, \dots,\bfx^{(N)}\}$.
 Here, were have $\bfy^{(n)} = (y_{11}^{(n)},\dots, y_{d^{(n)}1}^{(n)}, \dots,  y_{1q}^{(n)}, \dots, y_{d^{(n)}q}^{(n)} )$,  and $\bfx^{(n)} = \{\bfx^{(n)}_{11}, \dots,\bfx^{(n)}_{d^{(n)}1}, \dots,\bfx^{(n)}_{1q}, \dots,\bfx^{(n)}_{d^{(n)}q}\}$, $n = 1, \dots, N$. 
 Note that this sample can be thought of as observing $q$ different characteristics of $N$ groups of individuals, which can be of diverse sizes. Moreover, we denote by $\bfbeta$ the regression coefficients $(\bfbeta_{1}, \dots,\bfbeta_{q})$. 

Then, the complete loglikelihood (disregarding terms that do not depend on the parameters) is given by 
\begin{align}\label{eq:heirlog}
\begin{split}
	& l_c(\bfbeta,  \nu, \eta  ; \tilde{\bfy}, \tilde{\bfx} ) \\ &= \sum_{n = 1}^N \Bigg( \eta \log(\nu) - \log(\eta) + \eta \log(\theta_0^{(n)}) - \nu\theta_0^{(n)}  \\ 
	& \quad + \sum_{j = 1}^{q} \bigg(  \theta_0^{(n)} \log( \eta / \nu  ) - (\eta / \nu) \theta_j^{(n)} +  \sum_{i = 1}^{d^{(n)}} \Big( y_{ij}^{(n)} \bfx_{ij}^{(n)} \bfbeta_j - \theta_{j}^{(n)} E_{ij}^{{(n)}}e^{\bfx_{ij}^{(n)} \bfbeta_j} \Big)\bigg)\Bigg),
\end{split}
\end{align}
where we have denoted by $\theta^{(n)}_j$ the (unobserved) values of $\Theta_j$, $j = 0,1,\dots, d$.

\textbf{E-step.} The first step is to compute the conditional expectation of the loglikelihood given the observed data. To this end, let us denote by $\bfbeta^{(k)}$,  $\eta^{(k)}$, and $\nu^{(k)}$ the current parameters after $k$ iterations. Then, for the $(k + 1)$'th iteration, the problem reduces to computing $\E^{(k+1)} (\Theta_0 \mid \tilde{\bfy})$, $\E ^{(k+1)}(\log(\Theta_0 ) \mid \tilde{\bfy})$, and $\E^{(k+1)} (\Theta_j \mid \tilde{\bfy})$, $j = 1, \dots, q$, where $\E ^{(k+1)}(\cdot)$ denotes the expectation computed using the current parameters. 
 We perform the calculations for one (generic) data point $\bfy= (\bfy_1, \dots,  \bfy_q)$, given that for a sample of size $N$, we simply require to sum over all data points on those expressions. First, we note that  \eqref{eq:hiermar} yields
 \begin{align*}
   \E^{(k+1)} (\Theta_j \mid {\bfy}) = \frac{y_{\bullet j} + \E^{(k+1)} (\Theta_0 \mid {\bfy}) }{e_{\bullet j} + ({\eta^{(k)}}/{\nu^{(k)}})} .
  \end{align*}
 
The other conditional expectations follow directly from \eqref{eq:hiermix} and are given by
  \begin{align*}
   \E^{(k+1)} (\Theta_0 \mid {\bfy}) = \sum_m \bar{w}_m^{(k)} \frac{m + \eta^{(k)}}{\sum_{j = 1}^{q}\log\big(1 + (\nu^{(k)} / \eta^{(k)}) e_{\bullet j}\big) + \nu^{(k)} }
  \end{align*}
as well as 
 \begin{align*}
   &\E^{(k+1)} (\log(\Theta_0) \mid {\bfy})  \\ &= \sum_m \bar{w}_m^{(k)} \bigg( \Psi(m + \eta^{(k)}) - \log\Big(\sum_{j = 1}^{q}\log\big(1 + (\nu^{(k)} / \eta^{(k)}) e_{\bullet j}\big) + \nu^{(k)}\Big) \bigg) ,
  \end{align*}
where $\Psi(\cdot)$ denotes the digamma function. 

\textbf{CM-step.} With the previous conditional expectations at hand, we now proceed to maximize the conditionally expected loglikelihood with respect to the model's parameters. This will provide us with updated parameters $\bfbeta^{(k + 1)}$,  $\eta^{(k + 1)}$, and $\nu^{(k + 1)}$.
We start by describing how to find $\bfbeta^{(k+1)}$, which is the value of $\bfbeta$ maximizing  the following expression
\begin{align*}
	\sum_{j=1}^{q}\sum_{n= 1}^{N}  \sum_{i=1}^{d^{(n)}} \Big( y_{ij}^{(n)} \bfx_{ij}^{(n)} \bfbeta_j - \E^{(k+1)} (\Theta_j \mid \bfY = \bfy^{(n)}) E_{ij}^{{(n)}}e^{\bfx_{ij}^{(n)} \bfbeta_j} \Big).
\end{align*}
This expression can be seen as the loglikelihood of $q$ independent (conventional) Poisson regressions with offset values $\E^{(k+1)} (\Theta_j \mid \bfY = \bfy^{(n)}) E_{ij}^{{(n)}}$. Hence, $\bfbeta^{(k+1)} $ can be found easily by fitting $q$ generalized linear models.
Now, regarding $\eta^{(k + 1)}$ and $\nu^{(k + 1)}$, we first take partial derivatives in \eqref{eq:heirlog} with respect to $\eta$ and $\nu$, obtaining
\begin{align*}
	\frac{\partial}{\partial \nu } l_c(\bfbeta,  \nu, \eta  ; \tilde{\bfy}, \tilde{\bfx} ) \propto \frac{1}{q}\Big( \frac{\eta}{\nu} - \frac{1}{N}\sum_{n=1}^N \theta_0^{(n)} \Big) - \frac{1}{\nu}\frac{1}{N}\sum_{n=1}^N \theta_0^{(n)} + \frac{\eta}{\nu^2}\frac{1}{N}\frac{1}{q}\sum_{n=1}^N \sum_{j=1}^q \theta_j^{(n)}
\end{align*}
as well as 
\begin{align*}
	\frac{\partial}{\partial \eta } l_c(\bfbeta,  \nu, \eta  ; \tilde{\bfy}, \tilde{\bfx} ) &\propto \frac{1}{q}\Big( \log(\nu) - \Psi(\eta) + \frac{1}{N}\sum_{n=1}^N \log( \theta_0^{(n)}) \Big) \\ &\quad+ \frac{1}{\nu}\frac{1}{N}\sum_{n=1}^N \theta_0^{(n)} - \frac{1}{\nu}\frac{1}{N}\frac{1}{q}\sum_{n=1}^N \sum_{j=1}^q \theta_j^{(n)} .
\end{align*}
By fixing $\eta^{(k)}$ in the first expression (making this a CM-step) and plugging in the previously computed conditional expectations, we can explicitly find $\nu^{(k+1)}$ as
\begin{align*}
	\nu^{(k+1)}=\frac{\eta^{(k)} /q-\bar{\Theta}_{0, k}+\sqrt{\left( \eta^{(k)} / q-\bar{\Theta}_{0, k}\right)^2+4 \eta^{(k)} \bar{\Theta}_{0, k} \bar{\Theta}_k/ q} }{2 \bar{\Theta}_{0, k} / q} ,
\end{align*}
where 
\begin{align*}
	& \bar{\Theta}_{0, k} = \frac{1}{N}\sum_{n = 1}^N \E^{(k+1)} (\Theta_0 \mid \bfY = {\bfy}^{(n)})   \,, \\
	& \bar{\Theta}_k = \frac{1}{N}\frac{1}{q}\sum_{n = 1}^N\sum_{j = 1}^q \E^{(k+1)} (\Theta_j \mid \bfY = {\bfy}^{(n)}) .
\end{align*}
To compute the update parameter $\eta^{(k + 1)}$, we are in principle required to find the solution to $\frac{\partial}{\partial \eta } l_c = 0$. This equation does not possess a closed-form solution, so the use of some numerical optimizer is necessitated. 
{
We propose performing a single Newton-Raphson step, which is both computationally efficient and sufficient to ensure convergence of the ECM algorithm \citep[confer with][]{meng1993maximum}. Thus, the update is given by
}
\begin{align*}
	\eta^{(k + 1)}=\eta^{(k)}-\frac{\left(\log \nu^{(k+1)}-\Psi\left(\eta^{(k)}\right)+ \bar{\Theta}_{0, k}^{[\log ]}\right)\! / q+ \bar{\Theta}_{0, k}/ \eta^{(k)}- \bar{\Theta}_k/\nu^{(k+1)}}{-\Psi^{\prime}\left(\eta^{(k)}\right)\! / q- \bar{\Theta}_{0, k}/\!\left(\eta^{(k)}\right)^2}
\end{align*}
with
\begin{align*}
	\bar{\Theta}_{0, k}^{[\log ]} = \frac{1}{N	}\sum_{n = 1}^N \E^{(k+1)} (\log(\Theta_0) \mid \bfY = {\bfy}^{(n)}).
\end{align*}

\subsection{Bivariate phase-type mixing} \label{subsec:ph}
As previously mentioned, a limitation of the hierarchical Gamma model is that it only allows for positive correlation of the margins of the mixing vector. To address this, we turn our attention to a more flexible model capable of describing a wide range of dependence structures, namely multivariate phase-type distributions. The use of phase-type distributions as random effects was first explored in \cite{yslas2021fitting} in the context of univariate and multivariate frailty models. Here, we closely follow the methodologies outlined in \cite{yslas2021fitting} to derive explicit formulas of several functionals and an estimation algorithm in the context of multivariate mixed Poisson regression. 
We start by recalling the definition of a (univariate) phase-type distribution.
Consider a time-homogeneous Markov jump process $ ( J_t )_{t \geq 0}$ on the state space $\{1, \dots, p, p+1\}$, where states $1,\dots,p$ are transient and state $p+1$ is absorbing. 
Then, $ ( J_t)_{t \geq 0}$ has an intensity matrix on the form
\begin{align*}
	\mat{\Lambda}= \left( \begin{array}{cc}
		\bfT &  \bft \\
		\0 & 0
	\end{array} \right)\!, 
\end{align*}
where $\bfT $ is a $p \times p$ sub-intensity matrix and $\bft$ is a $p$-dimensional column vector satisfying $\bft =- \bfT \, \1$ with $\1 $ the $p$-dimensional column vector of ones.
 Assume that the process starts in a transit state $k$ with probability $ \pi_{k} = \P(J_0 = k)$, $k = 1,\dots, p$, and that $\sum_{k=1}^p \pi_k = 1$, that is, the process cannot start in the absorbing state. Finally, let $\bfpi = (\pi_1 ,\dots,\pi_p )$ be the vector of initial probabilities. Then, we say that the time until absorption $
	\Theta = \inf \{ t \geq  0 : J_t = p+1 \}$
has a phase-type distribution with representation $(\bfpi,\bfT )$, and we write $\Theta \sim \mbox{PH}(\bfpi,\bfT )$. 

Several functionals of phase-type distributions possess closed-form formulas in terms of matrices. For instance, the density function $f_\Theta$ of $\Theta \sim  \mbox{PH}(\bfpi , \bfT )$ is given by 
\begin{align*}
 f_\Theta(\theta) =&\vect{\pi}e^{  \mat{T} \theta} \vect{t} , \quad \theta>0 ,
\end{align*}
where the exponential of a matrix $\mat{M}$ is defined by
\begin{align*}
	e^{\mat{M}} = \sum_{r= 0}^{\infty} \frac{\mat{M}^{r}}{r!}.
\end{align*}
For our present purposes, we are particularly interested in the Laplace transform of $\Theta$, which can be computed as
\begin{align*}
	\mathcal{L}_\Theta(u) = \bfpi (u \mat{I} - \bfT )^{-1} \bft, \quad u \geq 0,
\end{align*}
with $\mat{I}$ the identity matrix of appropriate dimension. 
Besides their mathematical tractability, phase-type distributions are particularly appealing as modeling tools due to their denseness within the set of distributions on the positive real line. {Notably, they also include several important distributional families as special cases, such as Erlang distributions ({that is} Gamma distributions with integer shape parameters), general Coxian laws, and mixtures thereof.} We refer to \cite{Bladt2017} for a comprehensive account of phase-type distributions. 

Note, in particular, that~\eqref{eq:mmpr} under the assumption of $\Theta \sim \mbox{PH}(\bfpi, \bfT)$ takes the explicit form
\begin{align*}
	f_{\bfY }(\bfy ) = \Big(\prod_{i=1}^d \left(E_ie^{\bfX_i \bfbeta} \right)^{y_i} \Big)\frac{(\sum_{i=1}^d y_i)!}{y_1!\cdots y_d!} \bfpi \Big( \Big(\sum_{i=1}^d E_ie^{\bfX_i \bfbeta} \Big)\mat{I} - \bfT \Big)^{-\sum_{i=1}^d y_i-1} \bft .
\end{align*}

Unfortunately, the posterior distribution is no longer phase-type but rather falls within the more general class of matrix exponential distributions, that is, distributions with rational Laplace transform; for further details, see~\cite{yslas2021fitting}. Nevertheless, the model retains a high degree of mathematical tractability and takes full advantage of the modeling capabilities of phase-type distributions, a point we demonstrate in the subsequent discussion of~\eqref{eq:mmpr_gen} with bivariate phase-type mixing.
{
\begin{remark}
We note that one can alternatively employ other phase-type extensions as mixing distributions. In particular, certain inhomogeneous phase-type specifications introduced in \cite{AlbrecherBladt2019} admit closed-form Laplace transforms and can be extended to the multivariate setting via the transformation approach in \cite{albrecher2020fitting}. While such models are natural candidates, we focus here on the conventional phase-type case, as it provides the most direct extension of the independent Gamma mixing considered in \cite{furrer2019experience}, and already offers sufficient flexibility for our application.
\end{remark}
}

We now consider a specific class of multivariate phase-type distributions to describe the mixing component in \eqref{eq:mmpr_gen}. Here, we focus on  the bivariate class of phase-type distributions for which the joint density is given by
\begin{align}\label{eq:bivPH}
	f_{\bf\Theta}(\theta_1, \theta_2) = \bfeta e^{{\bfT_{11} \theta_1}}\bfT_{12}e^{{\bfT_{22} \theta_2}}(-\bfT_{22}) \1 , \quad \theta_1, \theta_2 > 0, 
\end{align}

where $\bfeta$ is a vector of initial probabilities of dimension $p_1$, $\bfT_{11}$ and $\bfT_{22}$ are sub-intensity matrices of dimensions $p_1 \times p_1$ and $p_2 \times p_2$, respectively, and $\bfT_{12}$ is a $p_1 \times p_2$ matrix satisfying $\bfT_{11} \1  = -\bfT_{12} \1 $. The motivation for considering this class is that it is dense on the set of distributions on $\R_+^2$  \citep[see][]{albrecher2020fitting}, thus allowing for the modeling of any dependence structure. {Some instances of bivariate specifications that fall under this framework include the case of independent phase-type components (in particular, independent Erlangs) and multivariate Erlang distributions constructed via mixtures of independent Erlangs.}
Moreover, this class can be estimated via an EM algorithm  \citep[see][]{ahlstrom1999parametric}.
The more recently introduced class of multivariate phase-type distributions in \cite{bladt2023tractable}, named the mPH class, possesses similar properties, and it is, in fact, a subfamily of the distributions with joint density \eqref{eq:bivPH}.
{
 Extensions beyond the bivariate case are possible, either through the construction outlined in \cite{AlbrecherBladtBladt2021} or via the mPH class, although such generalizations come with additional computational challenges \citep[confer with][for a detailed discussion]{yslas2021fitting}. In the present work, however, we focus only on the bivariate 
}
class described by~\eqref{eq:bivPH} to keep the presentation concise
{
and in accordance with the present application. 
}

We first note that~\eqref{eq:bivPH} yields
\begin{align*}
	 \E( \exp({-u_1 \Theta_1 - u_2\Theta_2})) =\bfeta ( u_1 \mat{I} - \bfT_{11})^{-1}\bfT_{12}( u_2 \mat{I} - \bfT_{22})^{-1}(-\bfT_{22}) \1 , \quad u_1, u_2 \geq 0 .
\end{align*}

Then, when $\bfTheta$ follows the distribution described by \eqref{eq:bivPH}, the above formula leads to the following explicit expression of~\eqref{eq:mmpr_gen}:
\begin{align*}
	 f_{\bfY}(\bfy) & = \left(  \prod_{i=1}^d \frac{\ (e_{i1})^{y_{i1}}}{y_{i1}!} \frac{ (e_{i2})^{y_{i2}}}{y_{i2}!}  \right)  \left( y_{\bullet 1} ! y_{\bullet 2} !\right) \times \\
	 & \qquad \bfeta \Big( e_{\bullet 1}\mat{I} - \bfT_{11}\Big)^{-y_{\bullet 1}-1}\bfT_{12}\Big(e_{\bullet 2} \mat{I} - \bfT_{22}\Big)^{-y_{\bullet 2} - 1}(-\bfT_{22}) \1  .
\end{align*}

Furthermore, it is also easy to see that for any $k, s \in \N_{0}$, it holds that
\begin{align}\label{eq:crossmom} 
\begin{split}
	&\E (\Theta_1^{k} \Theta_2^{s} \mid \bfY = \bfy) \\
	& = \frac{\big(y_{\bullet 1} + k\big)! \big(y_{\bullet 2} + s\big)!}{y_{\bullet 1}! y_{\bullet 2}!} \frac{\bfeta \Big( e_{\bullet 1}\mat{I} - \bfT_{11}\Big)^{-y_{\bullet 1}- k-1}\bfT_{12}\Big( e_{\bullet 2} \mat{I} - \bfT_{22}\Big)^{-y_{\bullet 2} -s - 1}(-\bfT_{22}) \1}{\bfeta \Big( e_{\bullet 1}\mat{I} - \bfT_{11}\Big)^{-y_{\bullet 1}-1}\bfT_{12}\Big( e_{\bullet 2} \mat{I} - \bfT_{22}\Big)^{-y_{\bullet 2} - 1}(-\bfT_{22}) \1}.
\end{split}
\end{align}
This last formula is particularly relevant for the estimation algorithm described below.

\subsubsection{Estimation} 
Maximum likelihood estimation of this model can also be performed via an EM algorithm, which we derive next. We start by writing down the complete loglikelihood, omitting terms which do not depend on the parameters of interest:
\begin{align}\label{eq:comll}
\begin{split}
	 l_c(\bfbeta, \bfeta, \bfT_{11} , \bfT_{12},  \bfT_{22} ; \tilde{\bfy}, \tilde{\bfx} )
	  	&  \propto \sum_{n= 1}^{N} \Bigg( \sum_{j=1}^{2}\sum_{i=1}^{d^{(n)}}\bigg(  y_{ij}^{(n)} \bfx_{ij}^{(n)} \bfbeta_j - \theta_{j}^{(n)} E_{ij}^{{(n)}}e^{\bfx_{ij}^{(n)} \bfbeta_j} \bigg) \\
	 &  \quad \qquad \quad + \log( f_{\bfTheta}(\bftheta^{(n)} ; \bfeta, \bfT_{11} , \bfT_{12},  \bfT_{22} ) ) \Bigg).
\end{split}
\end{align}

\textbf{E-step.} We denote by $\bfbeta^{(k)}$,  $\bfeta^{(k)}$, $\bfT_{11}^{(k)}$, $\bfT_{12}^{(k)}$, and $\bfT_{22}^{(k)}$  the current parameters after $k$ iterations. Then, for the $(k + 1)$'th iteration, we are required to compute $\E^{(k+1)} (\Theta_1 \mid \tilde{\bfy})$, $\E ^{(k+1)}(\Theta_2 \mid \tilde{\bfy})$, and $\E^{(k+1)} (\log( f_{\bfTheta}(\bfTheta)) \mid \tilde{\bfy}) $. 
 As before, it is enough to perform the calculations for one (generic) data point $\bfy= (\bfy_1, \bfy_2)$. Note that \eqref{eq:crossmom} yields
\begin{align*}
	&\E^{(k+1)} (\Theta_1 \mid \bfY = \bfy) \\
	&= \big(y_{\bullet 1} + 1\big)  \frac{\bfeta^{(k)} \Big( e_{\bullet 1}\mat{I} - \bfT_{11}^{(k)}\Big)^{-y_{\bullet 1}- 2}\bfT_{12}^{(k)}\Big( e_{\bullet 2} \mat{I} - \bfT_{22}^{(k)}\Big)^{-y_{\bullet 2} - 1}(-\bfT_{22}^{(k)}) \1}{\bfeta^{(k)} \Big( e_{\bullet 1}\mat{I} - \bfT_{11}^{(k)}\Big)^{-y_{\bullet 1}-1}\bfT_{12}^{(k)}\Big( e_{\bullet 2} \mat{I} - \bfT_{22}^{(k)}\Big)^{-y_{\bullet 2} - 1}(-\bfT_{22}^{(k)}) \1}
\end{align*}
as well as
\begin{align*}
	&\E^{(k+1)} (\Theta_2 \mid \bfY = \bfy) \\ 
	&= \big(y_{\bullet 2} + 1\big)  \frac{\bfeta^{(k)} \Big( e_{\bullet 1}\mat{I} - \bfT_{11}^{(k)}\Big)^{-y_{\bullet 1}- 1}\bfT_{12}^{(k)}\Big( e_{\bullet 2} \mat{I} - \bfT_{22}^{(k)}\Big)^{-y_{\bullet 2} - 2}(-\bfT_{22}^{(k)}) \1}{\bfeta^{(k)} \Big( e_{\bullet 1}\mat{I} - \bfT_{11}^{(k)}\Big)^{-y_{\bullet 1}-1}\bfT_{12}^{(k)}\Big( e_{\bullet 2} \mat{I} - \bfT_{22}^{(k)}\Big)^{-y_{\bullet 2} - 1}(-\bfT_{22}^{(k)}) \1}.
\end{align*}
For the last conditional expectation, we have in full generality that 
\begin{align*}
	&\E^{(k+1)}\big(\log( f_{\bfTheta}(\bfTheta; \bfeta, \bfT_{11} , \bfT_{12},  \bfT_{22})\big) \mid \bfY = \bfy) \\
	&= \int_{\R_+^2}   \log( f_{\bfTheta}(\bftheta ; \bfeta, \bfT_{11} , \bfT_{12},  \bfT_{22})) f_{\bfTheta | \bfY}(\bftheta | \bfy) \,\mathrm{d}\bftheta ,
\end{align*} 
where 
\begin{align*}	
f_{\bfTheta|\bfY}(\bftheta|\bfy ) =  \frac{f_{\bfY|\bfTheta}(\bfy|\bftheta )f_{\bfTheta}(\bftheta )}{f_{\bfY}(\bfy )}
	= \frac{\Big(\prod_{j = 1 }^{2} \frac{(\theta_j )^{y_{\bullet j} }}{y_{\bullet j} !}\exp(-\theta_j e_{\bullet j}  ) \Big) \bfeta e^{{\bfT_{11} \theta_1}}\bfT_{12}e^{{\bfT_{22} \theta_2}}(-\bfT_{22}) \1 }{\bfeta \Big( e_{\bullet 1} \mat{I} - \bfT_{11}\Big)^{-y_{\bullet 1} -1}\bfT_{12}\Big( e_{\bullet 2}  \mat{I} - \bfT_{22}\Big)^{-y_{\bullet 2} - 1}(-\bfT_{22}) \1} .
\end{align*}
A general closed-form solution to the integral above is not available. 
{
However, an explicit formula is not required, as the subsequent maximization of this expression with respect to the phase-type parameters can be translated into fitting a bivariate phase-type distribution to a known target distribution, as earlier noted in \cite{yslas2021fitting}. This approach eliminates the need for a closed-form solution, allowing the maximization to be carried out numerically, as discussed next.
}

\textbf{M-step.}
Note that, given the form of the complete loglikelihood~\eqref{eq:comll}, we can split the maximization into maximization of the regression parameter and maximization of the bivariate phase-type parameters, respectively. Regarding $\bfbeta^{(k+1)}$, it is clear from~\eqref{eq:comll} that the updated value can be found by performing two (conventional) Poisson regressions. On the other hand, to obtain the updated parameters $\bfeta^{(k+1)}$, $\bfT_{11}^{(k+1)}$, $\bfT_{12}^{(k+1)}$, and $\bfT_{22}^{(k+1)}$, we note that maximizing
\begin{align*}
\E^{(k+1)} (\log( f_{\bfTheta}(\bfTheta; \bfeta, \bfT_{11} , \bfT_{12},  \bfT_{22})) \mid \tilde{ \bfy})
\end{align*}
is equivalent to maximizing
{
\begin{align*}
& \frac{1}{N}  \E^{(k+1)} (\log( f_{\bfTheta}(\bfTheta; \bfeta, \bfT_{11} , \bfT_{12},  \bfT_{22})) \mid \tilde{ \bfy}) \\
&=  \int_{\R_+^2}   \log( f_{\bfTheta}(\bftheta ; \bfeta, \bfT_{11} , \bfT_{12},  \bfT_{22})) \left[ \frac{1}{N} \sum_{n = 1}^Nf_{\bfTheta | \bfY}(\bftheta | \bfy^{(n)}) \right] \, \mathrm{d}\bftheta
\end{align*} 
}
and that 
\begin{align}\label{eq:dens_op}	
\frac{1}{N} \sum_{n = 1}^N f_{\bfTheta|\bfY}(\bftheta|\bfy^{(n)} )
\end{align}
is a joint density function.
{
Finally, note that maximizing the above expression with respect to the phase-type parameters is equivalent to minimizing the Kullback--Leibler divergence between the bivariate phase-type distribution and \eqref{eq:dens_op}. The latter, in turn, can be carried out by fitting a bivariate phase-type model to the target distribution specified by \eqref{eq:dens_op}.
}
In~\cite{albrecher2020fitting}, it is shown how the EM algorithm in~\cite{ahlstrom1999parametric} can be modified for such a purpose. We also refer to~\cite{yslas2021fitting} for further details. 

\subsubsection{Numerical considerations}

When working with a large volume of data points, quantities such as $y_{\bullet j} = \sum_{i=1}^d y_{ij}$ can take considerably large values, 
{ 
and the direct evaluation of the previous expressions, which involve factorials, powers, and exponentials, may lead to combinations of extremely large and small numbers,
}
making the computation of some of the previous expressions unfeasible in the form given. However, this can be easily solved by working with alternative and equivalent expressions,
{ 
where such terms combine to result in numerically stable quantities.
}
More specifically, for the joint density, it is convenient to work with
\begin{align*}
	 f_{\bfY}(\bfy) & = \left(  \prod_{j=1}^2 \left( \prod_{i=1}^d  \frac{( e_{i j} )^{y_{ij}}}{y_{ij}!} \right) \frac{ y_{\bullet j} !}{(e_{\bullet j}  )^{y_{\bullet j} + 1}}\right) \times  \\
	  & \quad \bfeta \Big( \mat{I} - \bfT_{11} (e_{\bullet 1} )^{-1}\Big)^{-y_{\bullet 1} -1}\bfT_{12}\Big( \mat{I} - \bfT_{22} (e_{\bullet 2}  )^{-1}\Big)^{-y_{\bullet 2}  - 1}(-\bfT_{22}) \1  .
\end{align*}
Similarly, to compute the conditional expectations $\E^{(k+1)} (\Theta_j \mid \bfY = \bfy)$, $j = 1, 2$, it is better to obtain them using the following general expression:
\begin{align*}
	\E (\Theta_1^{k} \Theta_2^{s} \mid \bfY = \bfy) & = \frac{(y_{\bullet 1} + k)! (y_{\bullet 2} + s)!}{y_{\bullet 1}! y_{\bullet 2}!} (e_{\bullet 1})^{-k}(e_{\bullet 2})^{-s}\times  \\
	& \quad \frac{\bfeta \Big( \mat{I} - \bfT_{11} (e_{\bullet 1})^{-1}\Big)^{-y_{\bullet 1}- k-1}\bfT_{12} \Big( \mat{I} - \bfT_{22} (e_{\bullet 2})^{-1}\Big)^{-y_{\bullet 2} -s - 1}(-\bfT_{22}) \1}{\bfeta \Big( \mat{I} - \bfT_{11}(e_{\bullet 1})^{-1}\Big)^{-y_{\bullet 1}-1}\bfT_{12}\Big(  \mat{I} - \bfT_{22}(e_{\bullet 2})^{-1}\Big)^{-y_{\bullet 2} - 1}(-\bfT_{22}) \1}. 
\end{align*}
Finally, regarding $\E^{(k+1)} (\log( f_{\bfTheta}(\Theta_1, \Theta_2 ; \bfeta, \bfT_{11} , \bfT_{12},  \bfT_{22})) \mid \bfY = \bfy)$, it is convenient to rewrite $f_{\bfTheta|\bfY}(\bftheta|\bfy )$ as
 \begin{align*}	
	f_{\bfTheta|\bfY}(\bftheta|\bfy ) = \frac{\Big(\prod_{j = 1 }^{2}e_{\bullet j} \exp\!\big( y_{\bullet j}\log(\theta_j e_{\bullet j})-\log(y_{\bullet j}!) -\theta_j e_{\bullet j} \big) \Big)\bfeta e^{{\bfT_{11} \theta_1}}\bfT_{12}e^{{\bfT_{22} \theta_2}}(-\bfT_{22}) \1 }{\bfeta \Big( \mat{I} - \bfT_{11} (e_{\bullet 1})^{-1}\Big)^{-y_{\bullet 1}-1}\bfT_{12}\Big(  \mat{I} - \bfT_{22} (e_{\bullet 2} )^{-1}\Big)^{-y_{\bullet 2} - 1}(-\bfT_{22}) \1} .
\end{align*}

\section{Numerical study}\label{sec:numerical}

We now use simulated data to provide a detailed example of the practical implementation of the proposed methodologies for experience rating for disability insurance. We focus on the disability model depicted in Figure~\ref{fig:multi_state} with group effects. In particular, we delve into the shrinkage estimation of the disability and recovery rates. This hands-on approach will not only demonstrate the feasibility of the techniques proposed in the previous sections, but also provide a concrete understanding of their performance. All simulations and numerical computations were performed using the R programming language{, with the corresponding code available through the GitHub repository 
\url{https://github.com/jorgeyslas/bivph-experience-rating.git}.}

The numerical study is divided into three parts. In Subsection~\ref{subsec:sim}, we describe the methodology and approach taken to simulate an insurance portfolio with different groups, providing insights into the underlying assumptions and scenarios. Subsection~\ref{subsec:fit} focuses on the outcomes of the estimation of the models under study, presenting a detailed examination of the results and drawing comparisons between them. Finally, in Subsection~\ref{subsec:pred}, we explore the predictive performance of the models by examining the results in various iterations of the simulation study and assessing their effectiveness for pricing. {The insights from the last subsection are particularly important from a pratical point of view.}

\subsection{Simulation} \label{subsec:sim}

\begin{figure}[!htbp]
\centering
\includegraphics[width=0.8\textwidth]{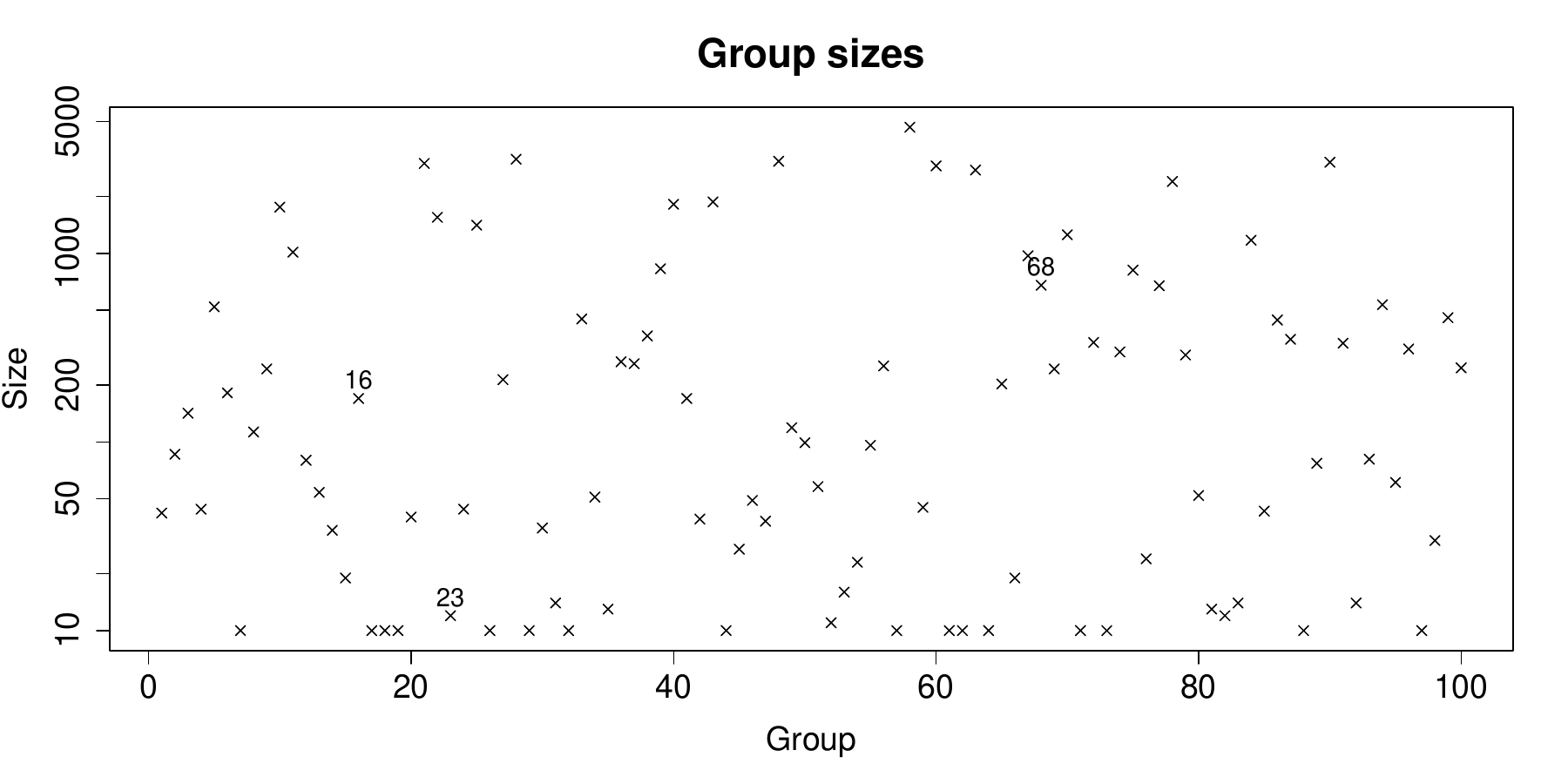}
\caption{Group sizes; groups $16$, $23$, and $68$ have been highlighted.}
\label{fig:group_sizes}
\end{figure}

We begin our numerical study by simulating an insurance portfolio comprised of $100$ distinct groups and a total of $50{,}000$ insured, each having a coverage period of $3$ years{, but also a retirement age (maximal contract time) of $67$}. Recognizing that real-world portfolios often comprise groups of varying sizes, we have {-- partly by hand --} allocated the insured to groups of differing sizes, confer with Figure~\ref{fig:group_sizes}. The group sizes {seem to follow}, approximately, a Gamma distribution with shape parameter {of about} $0.05$ and rate parameter {of about} $0.0001$.

\begin{figure}[!htbp]
\centering
\includegraphics[width=0.8\textwidth]{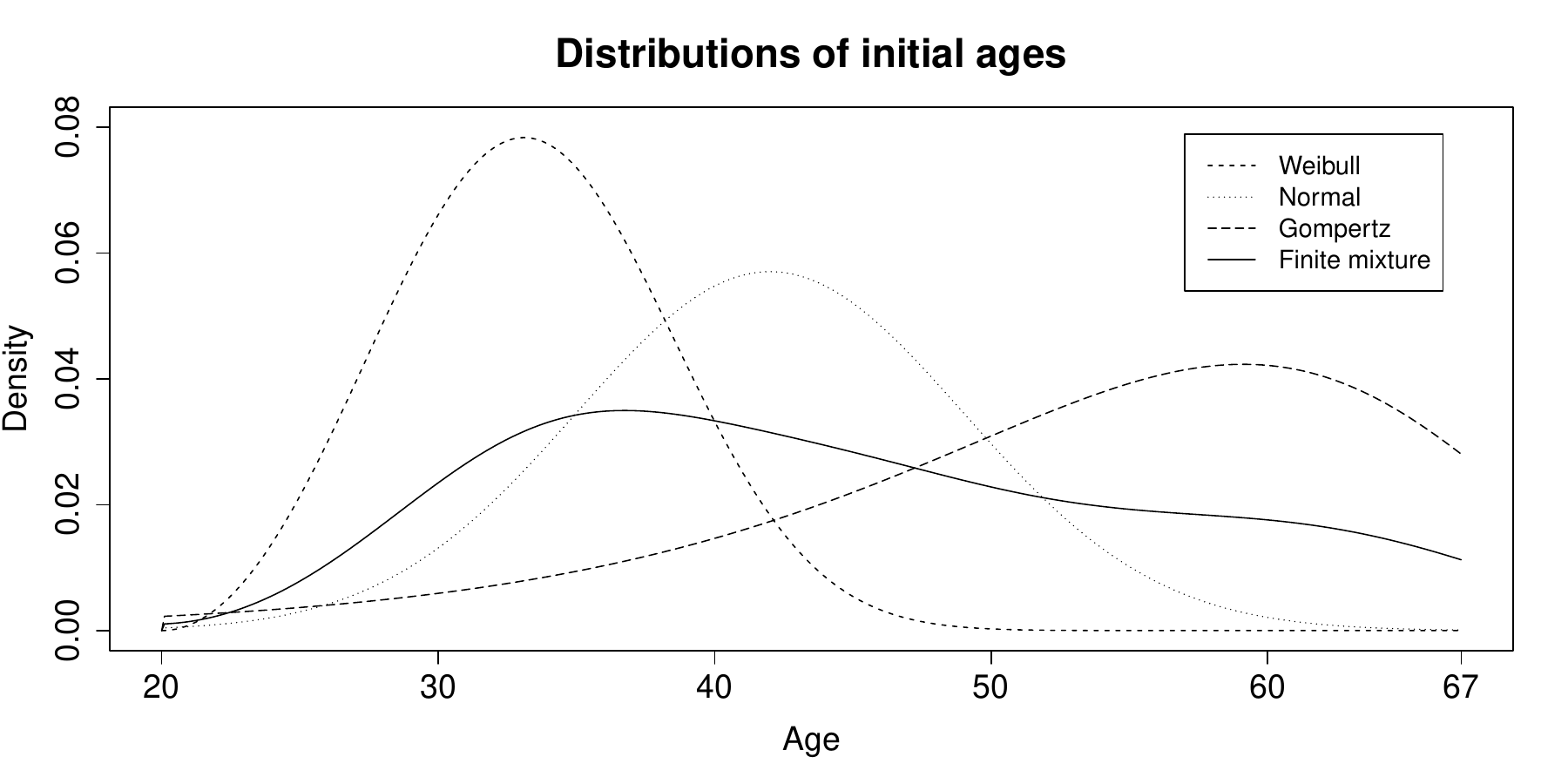}
\includegraphics[width=0.39\textwidth]{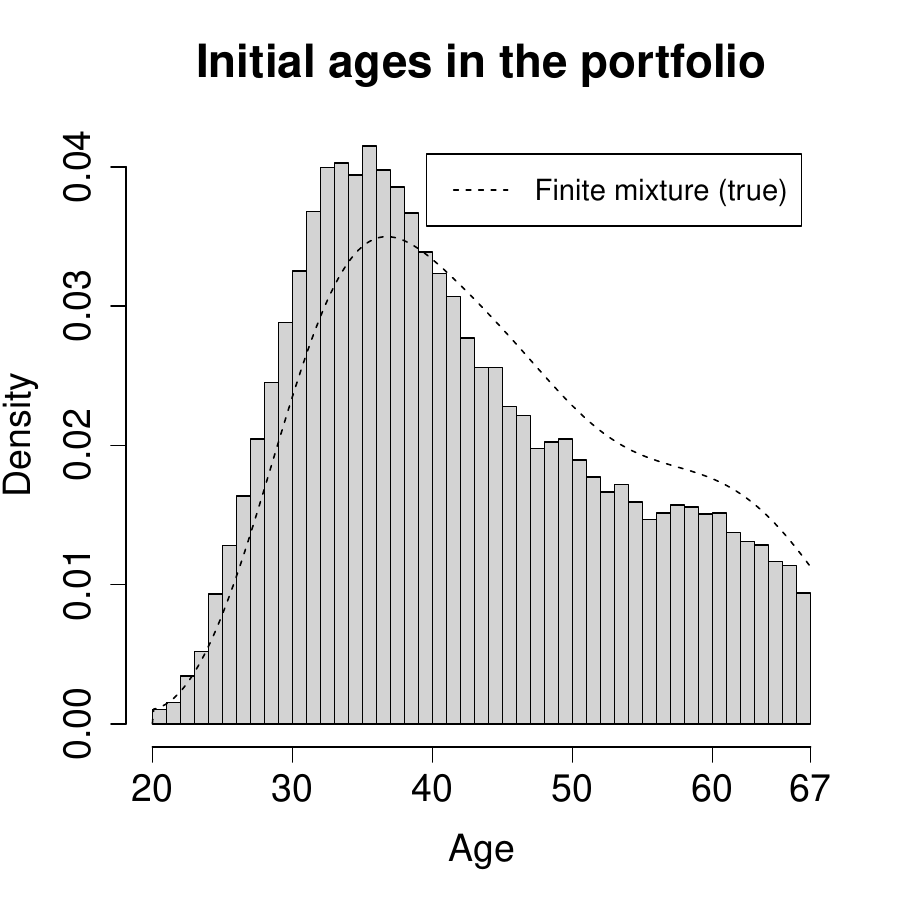}
\includegraphics[width=0.39\textwidth]{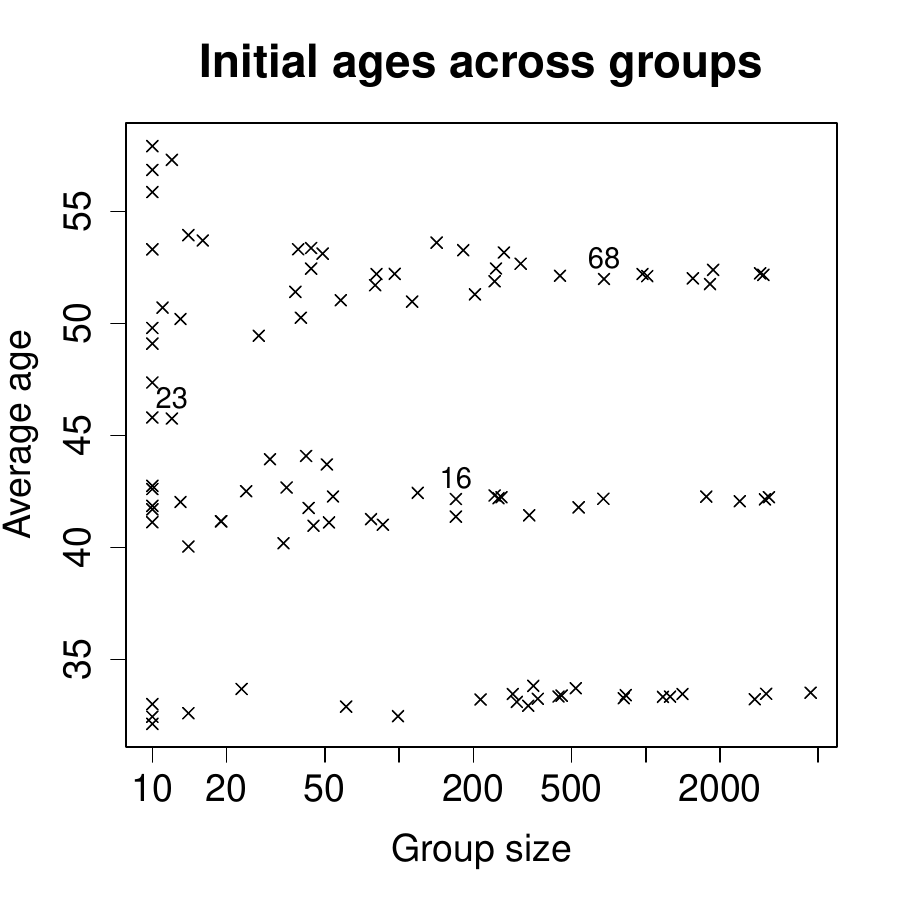}
\caption{Mixture density for the initial ages (top), histogram of simulated ages against original density (bottom left), and average initial age per group size with groups $16$, $23$, and $68$ highlighted (bottom right). }
\label{fig:age_dists}
\end{figure}

Next, we focus on the initial age composition within each group. To model the initial ages of the different groups, we randomly generate entry ages between $20$ and $67$ years following a mixture distribution. The choice of a mixture distribution is to reflect different age profiles in the portfolio, while the age range is representative of an actual insurance portfolio. More specifically, for the mixture distribution, we consider the following three distinct mixture components: A Weibull distribution with shape parameter $3$ and scale parameter $15$, representing younger groups,  a Normal distribution with mean $47$ and standard deviation $7$, symbolizing middle-aged groups, and a Gompertz distribution with shape parameter $0.1$ and rate parameter $0.002$, characterizing older groups. The corresponding mixture weights are $0.25$, $0.35$, and $0.4$, respectively. {Ages are drawn from the mixture distribution conditional on the range from $20$ to $67$ years (via rejection and resampling).} Figure~\ref{fig:age_dists} depicts the described mixture distribution along with a visual representation of the simulated initial ages. 

We proceed with the specification of the transition rates. For the fixed mortality rates, we select the same intensitites as in~\cite{furrer2019experience}:
\begin{align*}
	& \mu_{ad}(t) = 0.0005 + 10^{5.88 + 0.038 t - 10} ,\\
	& \mu_{id}(t) = \exp(-7.25 + 0.07 t) .
\end{align*}

On the other hand, for the disability and recovery rate, we consider the following baseline rates 
\begin{align*}
	& \mu_{ai}(t) = \exp(-4.5 - 0.018  t + 0.00064  t^2) , \\
	& \mu_{ia}(t) =\exp(0.3 - 0.049 t) .
\end{align*}

\begin{figure}[!htbp]
\centering
\includegraphics[width=0.44\textwidth]{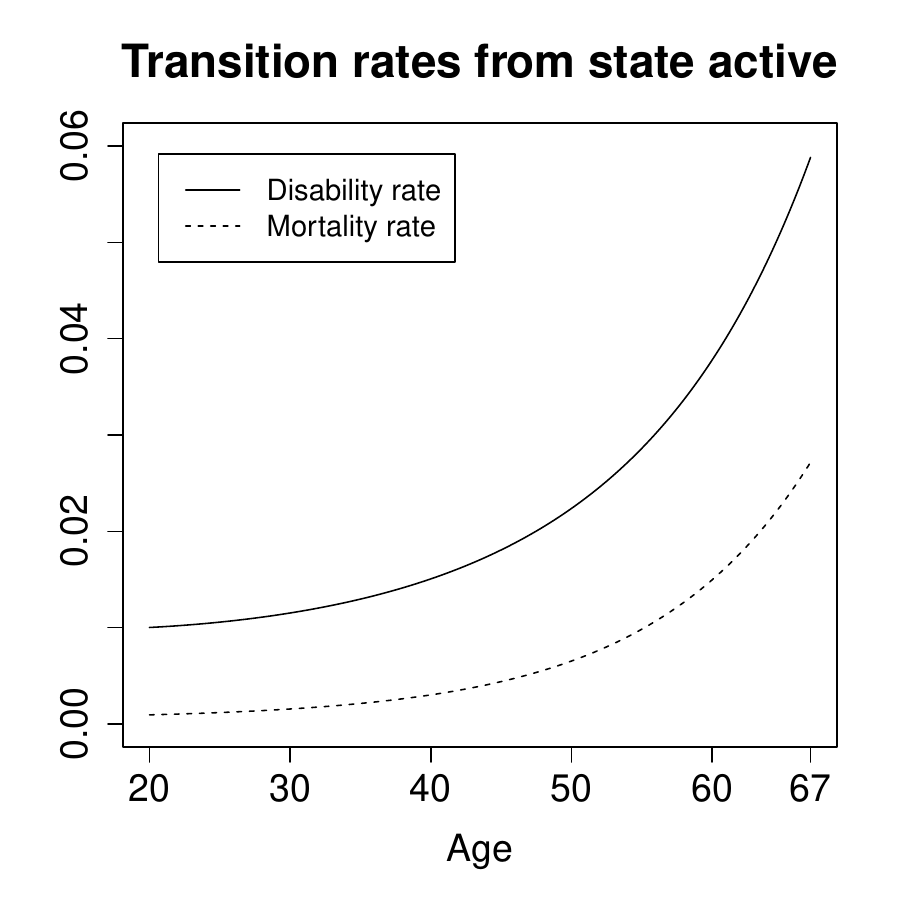}
\includegraphics[width=0.44\textwidth]{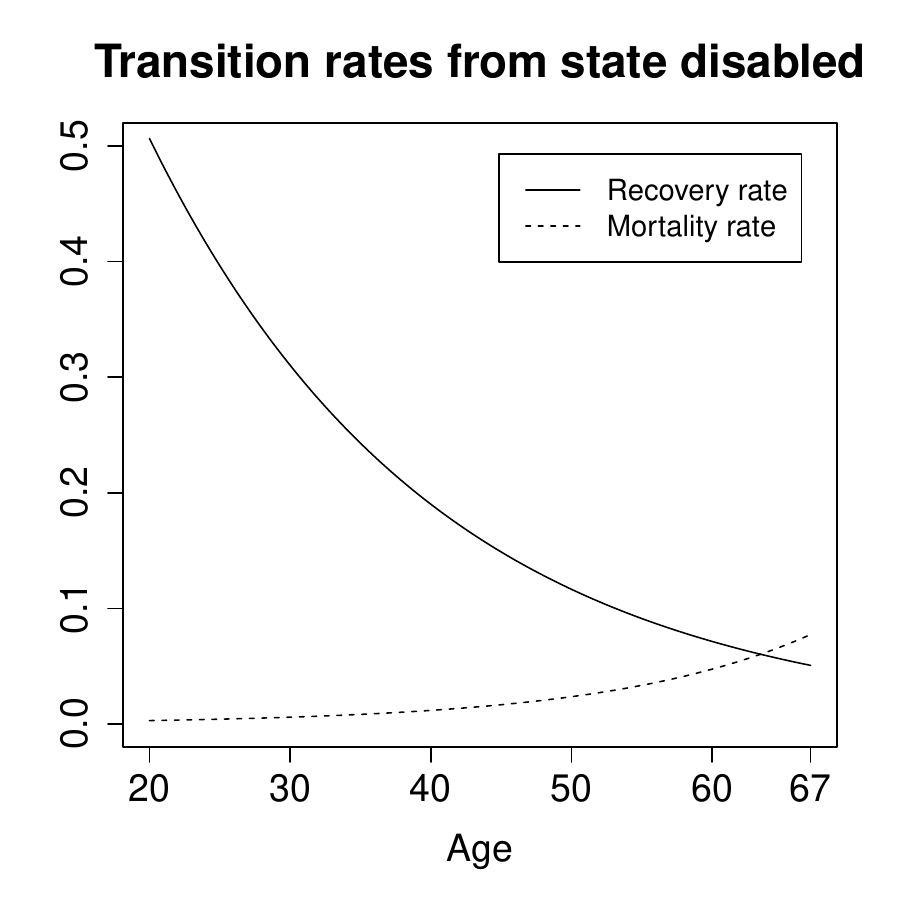}
\caption{True transition rates for the Markov chain.}
\label{fig:true_rates}
\end{figure}

Figure~\ref{fig:true_rates} shows the transition rates just defined. It is worth mentioning that we keep the above-described initial group sizes, age distributions, and baseline transition rates fixed throughout the study.

\begin{figure}[!htbp]
\centering
\includegraphics[width=0.9\textwidth]{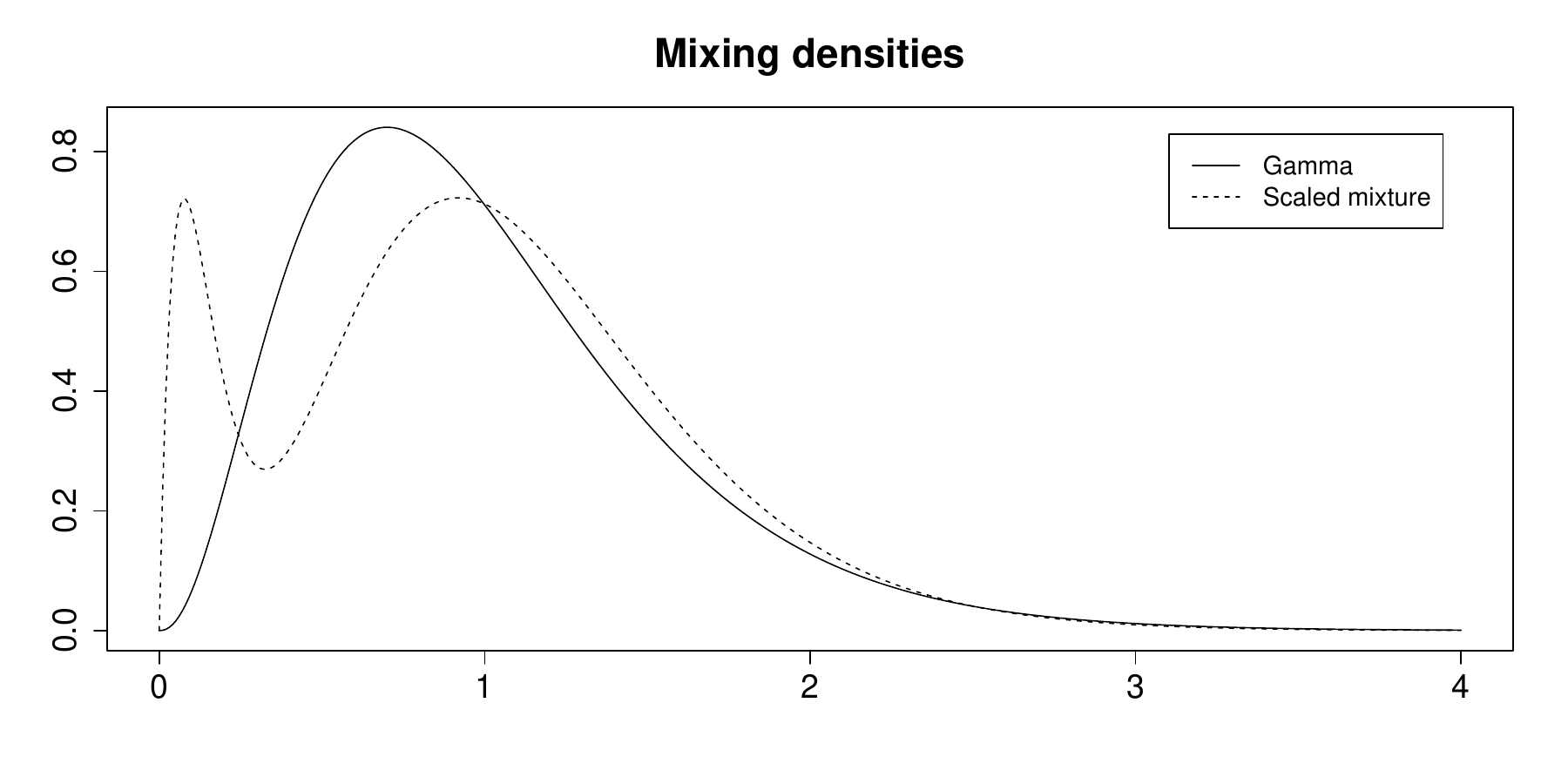}
\caption{Marginal densities of the mixing vector $\bfTheta$. }
\label{fig:priors}
\end{figure}

We next introduce a group effect $\bfTheta = (\Theta_{ai}, \Theta_{ia})$ with $\Theta_{ai}$ affecting the disability rate and $\Theta_{ia}$ affecting the recovery rate.  Specifically, we explore two distinct marginal distributions: a $\text{Gamma}(10/3, 10/3)$-distribution and a mixture of two Gamma distributions, namely $\text{Gamma}(5, 2)$ and $\text{Gamma}(2, 6)$ with weights $0.85$ and $0.15$, which is then scaled to have mean equal to one. The corresponding densities of these distributions are depicted in Figure~\ref{fig:priors}. 

One of the primary goals of this study is to investigate how the proposed models perform with respect to different dependence structures of the mixing vector $\bfTheta$. Thus, we consider three possible scenarios: independence, positive correlation, and negative correlation. This setting, along with the two distinct marginal specifications, enables us to evaluate the robustness and flexibility of the reviewed models in capturing complex joint behaviors, reflecting various possible real-world situations for disability insurance portfolios. More specifically, to introduce dependency within the random vector $\bfTheta$, we employ a Clayton copula with three different values of Kendall's tau, $\rho_{\tau}$: $-0.5$, $0$, and $0.5$, thus encompassing all three dependency cases.

To keep the presentation concise, the presentation focuses on the following specific scenaries for the joint prior distribution of $\bf\Theta$:
\begin{itemize} 
\item[] \textbf{Case A}: \textit{Independent Gamma }. Here, $\Theta_{ai}$ and $\Theta_{ia}$ are taken independent $\text{Gamma}(10/3, 10/3)$-distributed. \\[-2mm]
\item[] \textbf{Case B}: \textit{Positive Correlated Scaled Mixture}. Here, we use the scaled mixture for the marginal behavior of both $\Theta_{ai}$ and $\Theta_{ia}$, bundled  together by a Clayton copula with $\rho_{\tau} = 0.5$. 
\end{itemize} 
Although only the phase-type model is able to capture negative correlation, making this case the least compelling for a direct comparison, we also consider the following scenario:
\begin{itemize}
\item[] \textbf{Case C}: \textit{Negative Correlated Scaled Mixture}. Here, we use the scaled mixture for the marginal behavior of both $\Theta_{ai}$ and $\Theta_{ia}$, but bundled together by a Clayton copula with $\rho_{\tau} = -0.5$. 
\end{itemize}

{
Note that only Case A coincides with one of the considered models (independent Gamma). In Cases B and C, none of the models under comparison matches the true specification, in particular because the dependence structure is generated using a Clayton copula. This choice ensures that the models' performance is assessed outside their own families as to not introduce unwanted bias in favor of one of the models. However, it is essential to point out that we purposely chose the same distribution for both margins, as the hierarchical model is restricted to having equal marginal variances. As shown at the end of the study, this restriction represents a significant limitation for the predictive performance of the hierarchical specification.
}

The portfolio data is simulated as follows. First, for each case (A, B, and C), values of the group effects are simulated across groups. Next, for each insured in each group, we simulate the path of a Markov jump process with state active at entry age and with appropriate transition rates. To take into account the coverage period of $3$ years, simulations are terminated at time $\inf\{t \geq 3 : Z_t = a\}$. Simulation of a path involves the simulation of inhomogeneous Poisson processes, which we carry out using the acceptance-rejection method of~\cite{lewis1979simulation}. In Table~\ref{tbl:sim_theta}, we report the simulated group effects for groups $16$, $23$, and $68$. These groups have been selected to represent a range of sizes, as can be seen in Figure~\ref{fig:group_sizes}. Specifically, we have purposefully chosen group $16$ to represent a small size, group $23$ for a medium size, and group $68$ to represent a large size group. The bottom right panel of Figure~\ref{fig:age_dists} provides the average initial age for each of these groups, further emphasizing their distinct characteristics.

\begin{table}[!htbp]
\centering
\begin{tabular}{l |c c c | c c c}
\hline
\multirow{2}{*}{Group (g)} & \multicolumn{3}{c|}{$\Theta_{ai}^{(g)}$ (disability rate)} & \multicolumn{3}{c}{$\Theta_{ia}^{(g)}$ (recovery rate)} \\
\cline{2-7}
	& Case A   & Case B   & Case C & Case A   & Case B   & Case C \\
 \hline 
 $16$ & $0.821751$   & $0.517464$   & $0.175590$  & $0.732254$ & $1.043226$ & $1.214685$ \\
 $23$ & $0.642414$   & $0.076484$   & $0.944184$  & $0.363231$ & $0.306847$ & $0.666432$  \\
 $68$ & $2.318506$   & $0.167267$   & $2.095090$  & $0.979070$ & $0.405047$ & $0.936463$  \\
\hline
\end{tabular}
\caption{True (simulated) group effects for the three selected groups.}
\label{tbl:sim_theta}
\end{table}

{Throughout, we keep the simulated group effects fixed -- and the various models are compared in their ability to assess these simulated effects. Consequently, to understand the performance of the different models across scenarios, one needs to have in mind the empirical, rather than the underlying, distribution of these effects. Therefore, we report in Table~\ref{tbl:empirical_dist_group_effects} empirical estimates of the mean, standard deviation, and Kendall's $\rho_\tau$ across cases.

\begin{table}[!htbp]
\centering
\begin{tabular}{l |c c c | c c c}
\hline
\multirow{2}{*}{Empirical} & \multicolumn{3}{c|}{$\Theta_{ai}$ (disability rate)} & \multicolumn{3}{c}{$\Theta_{ia}$ (recovery rate)} \\
\cline{2-7}
	 & Case A   & Case B   & Case C & Case A   & Case B   & Case C \\
 \hline 
 Mean & $0.994437$   & $0.853114$   & $\phantom{-}0.988569$  & $1.003830$ & $0.916206$ & $\phantom{-}0.955082$ \\
 Std.\ deviation& $0.612241$   & $0.584568$   & $\phantom{-}0.552639$  & $0.594999$ & $0.512831$ & $\phantom{-}0.476886$  \\
 Kendall's $\rho_\tau$ & $0.061818$   & $0.565657$   & $-0.452525$  & $0.061818$ & $0.565657$ & $-0.452525$  \\ 
\hline 
\end{tabular}
\caption{Empirical mean, standard deviation, and Kendall's $\tau$ for the simulated group effects.}
\label{tbl:empirical_dist_group_effects}
\end{table}
}


\subsection{Estimation} \label{subsec:fit}

Having simulated a dataset for Case A, Case B, and Case C, respectively, we now proceed to perform shrinkage estimation of the disability and recovery rates using the independent Gamma, which we also call `simple' following~\cite{furrer2019experience}, the hierarchical Gamma, and the bivariate phase-type model. 
{
It is essential to note that we set $p_1 = p_2 = 3$ for the bivariate phase-type mixing in Case A, while for Cases B and C we use $p_1 = p_2 = 6$. These dimensions were chosen to provide sufficient flexibility to capture the different dependence structures considered -- particularly in Cases B and C -- while keeping the estimation computationally feasible and ensuring stability in the subsequent predictive performance study based on multiple portfolio realizations. For practical applications, however, one may follow the standard model selection strategy used in phase-type fitting, beginning with low-dimensional specifications and gradually increasing complexity, with improvements assessed through, for instance, loglikelihood comparisons, which can be additionally supported with visual diagnostic tools \citep[see][for a detailed discussion]{yslas2021fitting}.
}

Across all models, we assume the following parametric form of the baseline transition rates:
\begin{align*}
	& \mu_{ai}(t;\bfbeta^{(ai)} ) = \exp(\beta_{0}^{(ai)} + \beta_1^{(ai)} t + \beta_2^{(ai)}  t^2) , \\
	& \mu_{ia}(t; \bfbeta^{(ia)}) =\exp(\beta_0^{(ia)} + \beta_1^{(ia)} t ) ,
\end{align*}
where $\bfbeta^{(ai)} = (\beta_{0}^{(ai)}, \beta_1^{(ai)}, \beta_2^{(ai)} )$ and $\bfbeta^{(ia)} = (\beta_0^{(ia)}, \beta_1^{(ia)})$. Note that this mirrors the parametrization of the original transition rates. 

Recall that shrinkage estimation consists of first performing maximum likelihood estimation of the corresponding mixed Poisson regression models via the EM algorithms presented in Section~\ref{sec:mpr}, thus obtaining different estimators of the regression parameters, ${\bfbeta}^{(ai)}$ and ${\bfbeta}^{(ia)}$, and the parameters of the specifications of $\bfTheta$. 
{
As it is customary in any form of EM algorithm, the number of steps is determined such that the absolute change in the parameter estimates or loglikelihoods between successive iterations becomes negligible.
}
Then, the shrinkage estimates of $\bfTheta$ are obtained by computing the Bayes estimators under quadratic loss utilizing the previously computed maximum likelihood estimates. Finally, if we denote by $\hat{\bfTheta}^{(g)} = (\hat{\Theta}_{ai}^{(g)}, \hat{\Theta}_{ia}^{(g)})$ the resulting shrinkage estimates for each group, $g =1,\dots, 100$, and by $\hat{\bfbeta}^{(ai)}$ and $\hat{\bfbeta}^{(ia)}$ the estimates of the regression parameters, the `plug-in' approach uses as estimates of the disability and recovery rates exactly
\begin{align*}
	t \mapsto \hat{\Theta}_{ai}^{(g)} \mu_{ai}(t;\hat{\bfbeta}^{(ai)} ) \quad \text{and} \quad t \mapsto \hat{\Theta}_{ia}^{(g)} \mu_{ia}(t; \hat{\bfbeta}^{(ia)})
\end{align*}
for all $g = 1,\dots, 100$. 

For the sake of completeness, we also explore two additional models in line with~\cite{furrer2019experience}. The first is called \textit{standard}. This approach disregards group effects altogether, providing a baseline for comparison. The second is called \textit{fixed effect}. Here, the group effects are estimated as fixed effects, disregarding shrinkage completely. Estimation in these last two models is straightforward in R using the \textit{glm()} function.

\begin{figure}[!htbp]
\centering
\includegraphics[width=0.44\textwidth]{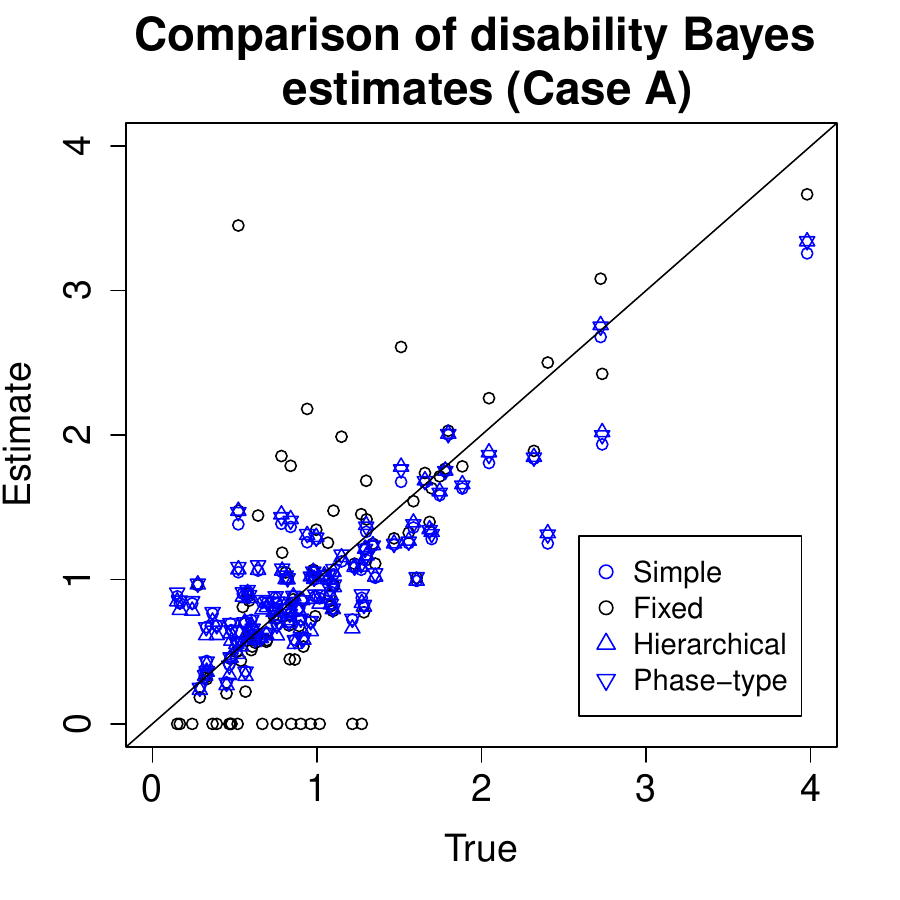}
\includegraphics[width=0.44\textwidth]{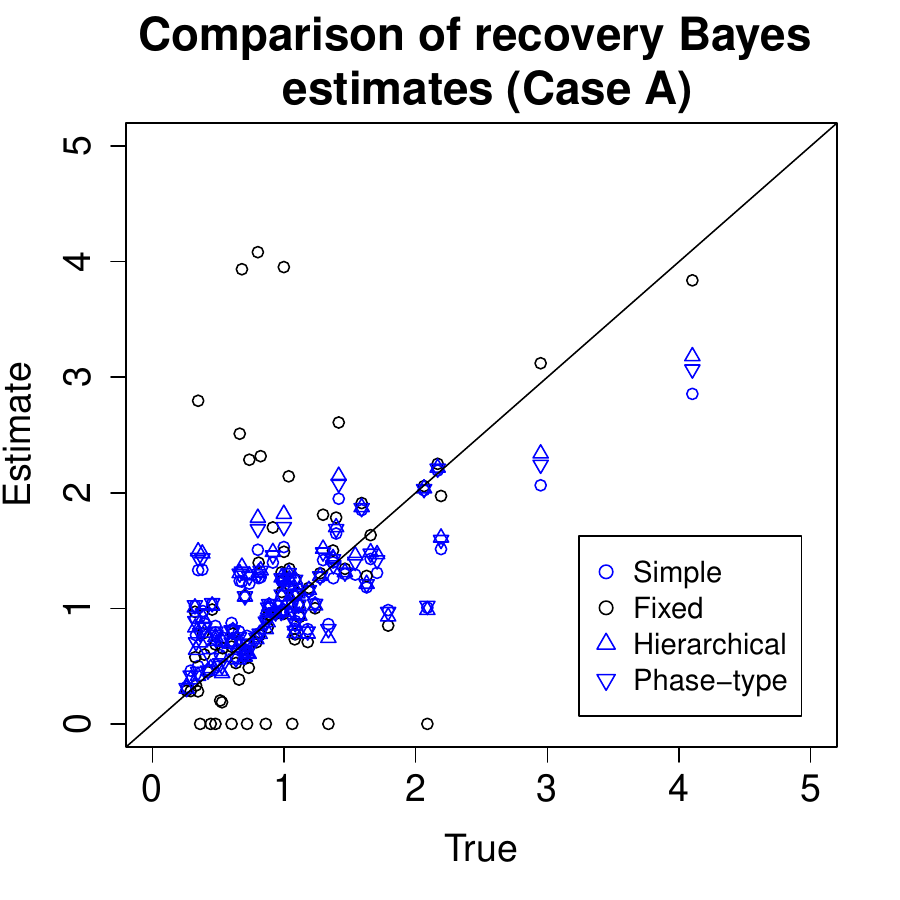}
\includegraphics[width=0.44\textwidth]{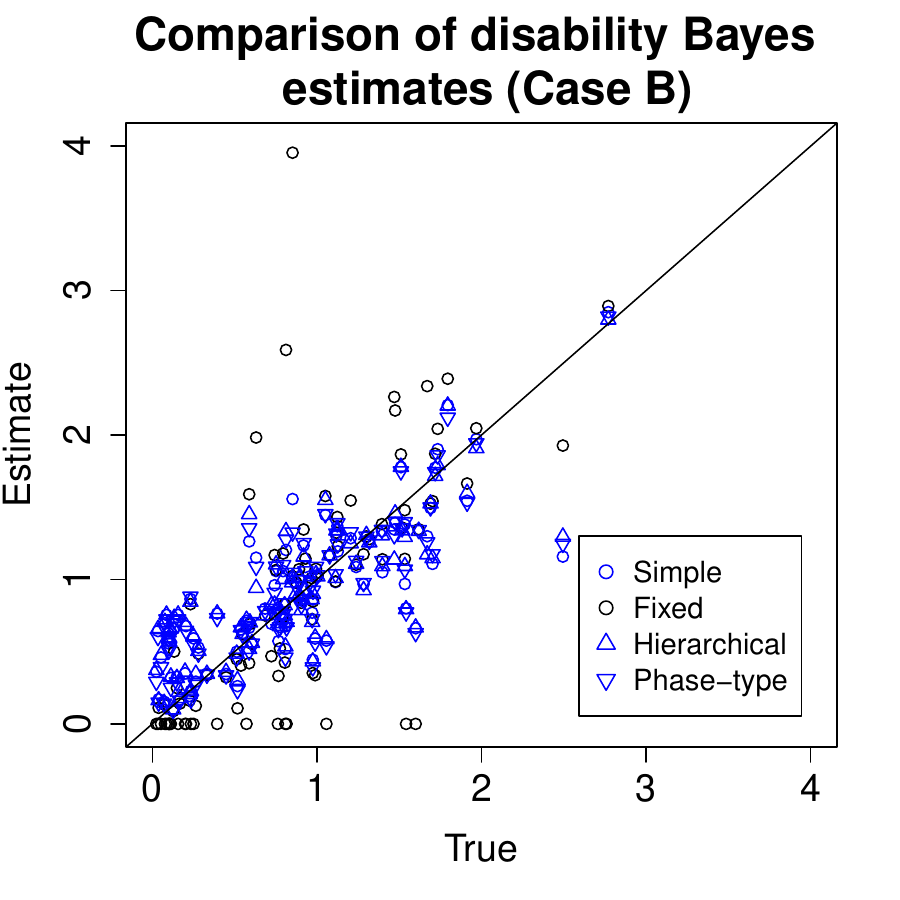}
\includegraphics[width=0.44\textwidth]{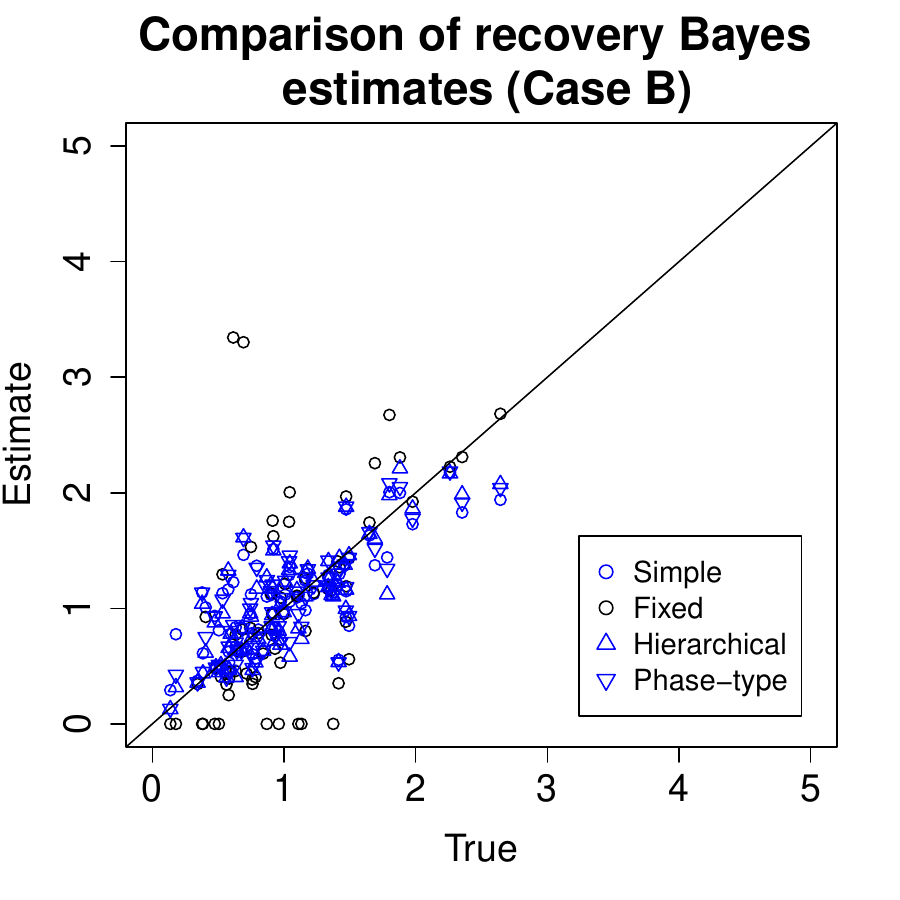}
\includegraphics[width=0.44\textwidth]{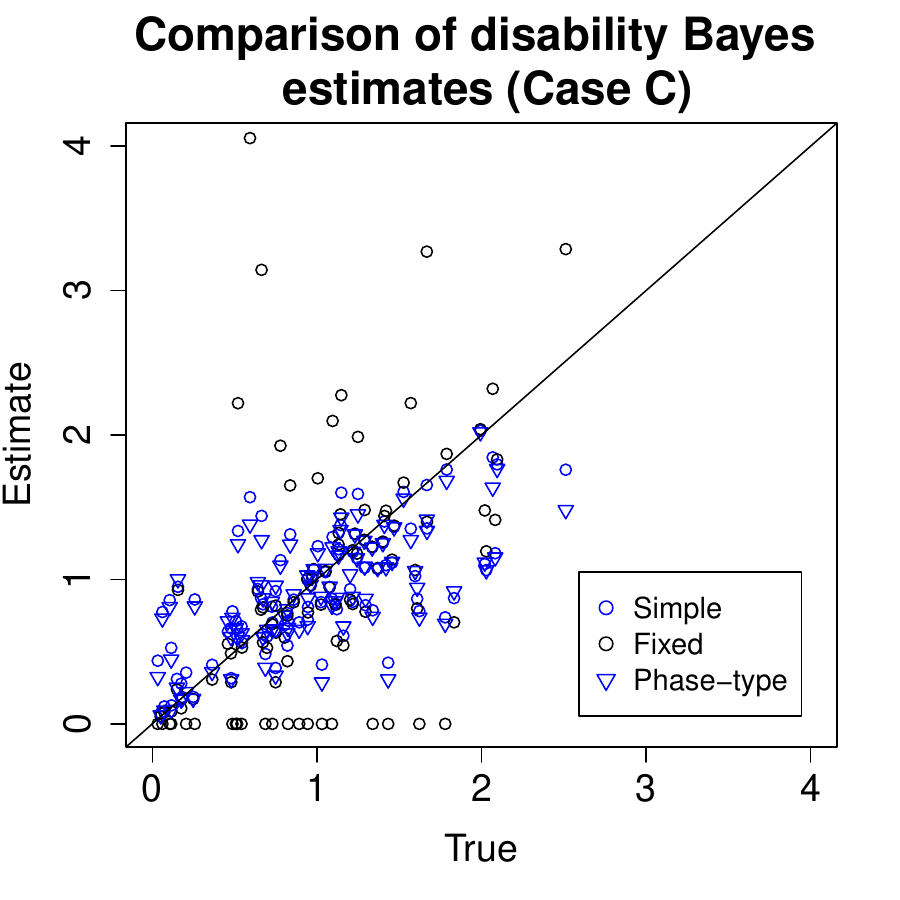}
\includegraphics[width=0.44\textwidth]{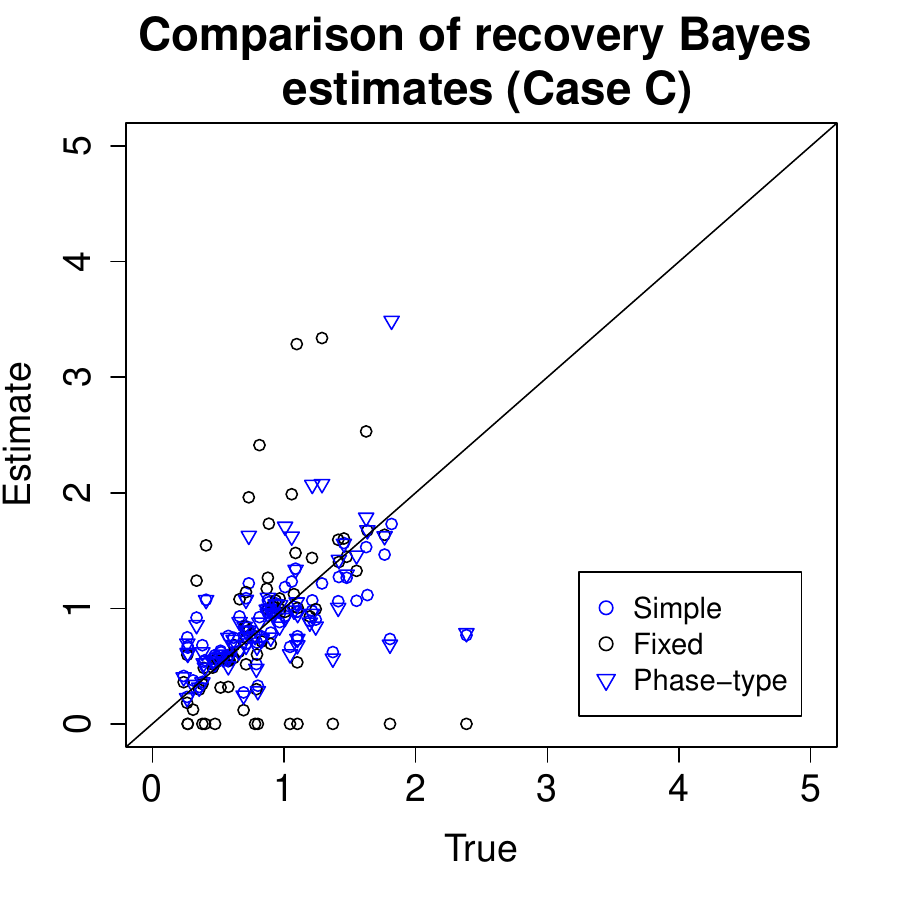}
\caption{Estimated group effects against true (simulated) values.}
\label{fig:theta_all}
\end{figure}

For a visual comparison of the fitted models, we would like to begin by providing plots of the original (simulated) group effects against the estimated ones. Nevertheless, given the presence of an intercept in the regression and also due to the well-known identifiability issue of phase-type distributions \citep[confer with][]{yslas2021fitting}, it becomes more suitable to graph the original group effects against the following scaled estimates:
\begin{align*}
	\hat{\Theta}_{ai}^{(g)} \frac{\sum_{k = 1}^{K} E_a^{k, (g)} \mu_{ai}(t_k;\hat{\bfbeta}^{(ai)} )}{\sum_{k = 1}^{K} E_a^{k, (g)} \mu_{ai}(t_k;{\bfbeta}^{(ai)}_{\text{true}} )}  \quad \text{and} \quad \hat{\Theta}_{ia}^{(g)} \frac{\sum_{k = 1}^{K} E_i^{k, (g)} \mu_{ia}(t_k;\hat{\bfbeta}^{(ia)} )}{\sum_{k = 1}^{K} E_i^{k, (g)} \mu_{ia}(t_k;{\bfbeta}_{\text{true}}^{(ia)} )}, 
\end{align*}
where ${\bfbeta}_{\text{true}}^{(ai)} = (-4.5, -0.018, 0.00064 )$, ${\bfbeta}_{\text{true}}^{(ai)}= (0.3, -0.049)$, and $E_a^{k, (g)}$ and $E_i^{k, (g)}$ are the relevant exposures across groups $g =1, \dots, 100$.  This procedure allows for a more meaningful comparison as it mitigates the aforementioned effects. The resulting plots for the group effects of the disability and recovery rates
{ across Cases A, B, and C are collected in Figure~\ref{fig:theta_all}, where we have omitted the hierarchical model in Case C, since it only allows for positive correlation and therefore cannot be meaningfully compared in a setting with negative dependence. For this reason, the hierarchical specification is consistently excluded from the Case C analysis throughout the study. } 
We observe that the models specifying a random effect tend to outperform the fixed effect one. This performance difference is particularly noticeable in the cases of groups without observed disabilities, leading to an estimation of a group effect equal to zero in the fixed effect case, while this is not the case for the random effects models. However, from the plots alone, it is not evident which model for $\bfTheta$ performs best in each case. Therefore, we now turn our attention to the analysis of the estimates of the transition rates. We focus on the groups $16$, $23$, and $68$.

%
%

\begin{figure}[!htbp]
\centering
\includegraphics[width=0.44\textwidth]{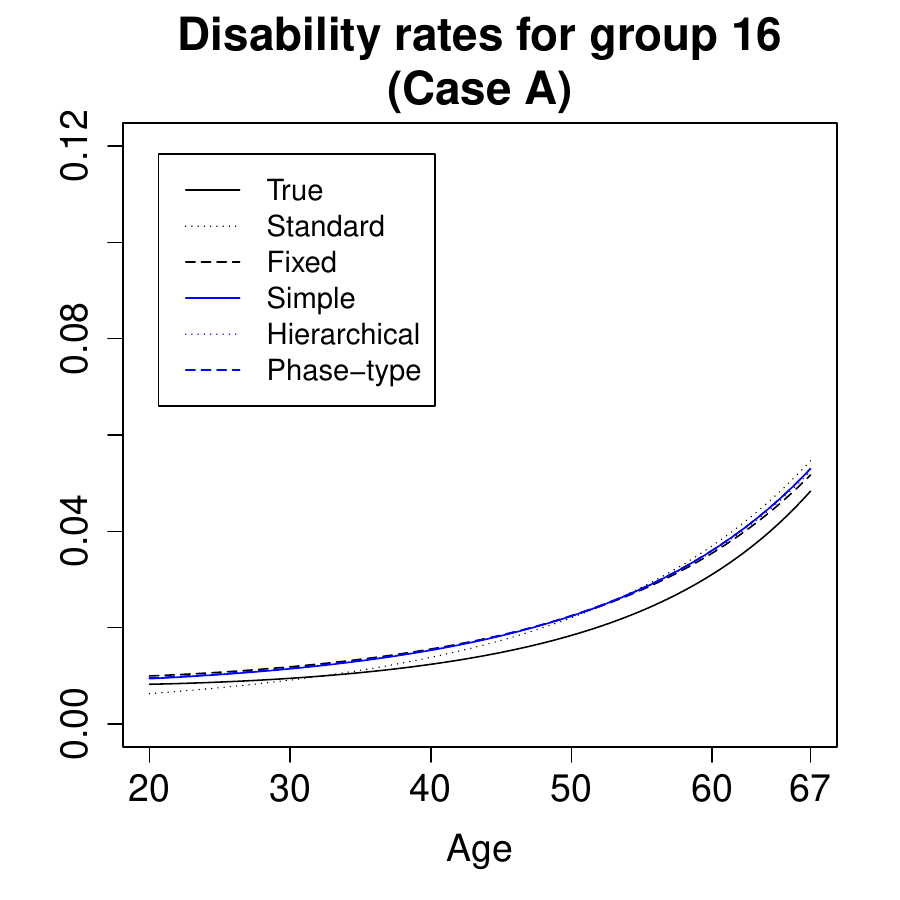}
\includegraphics[width=0.44\textwidth]{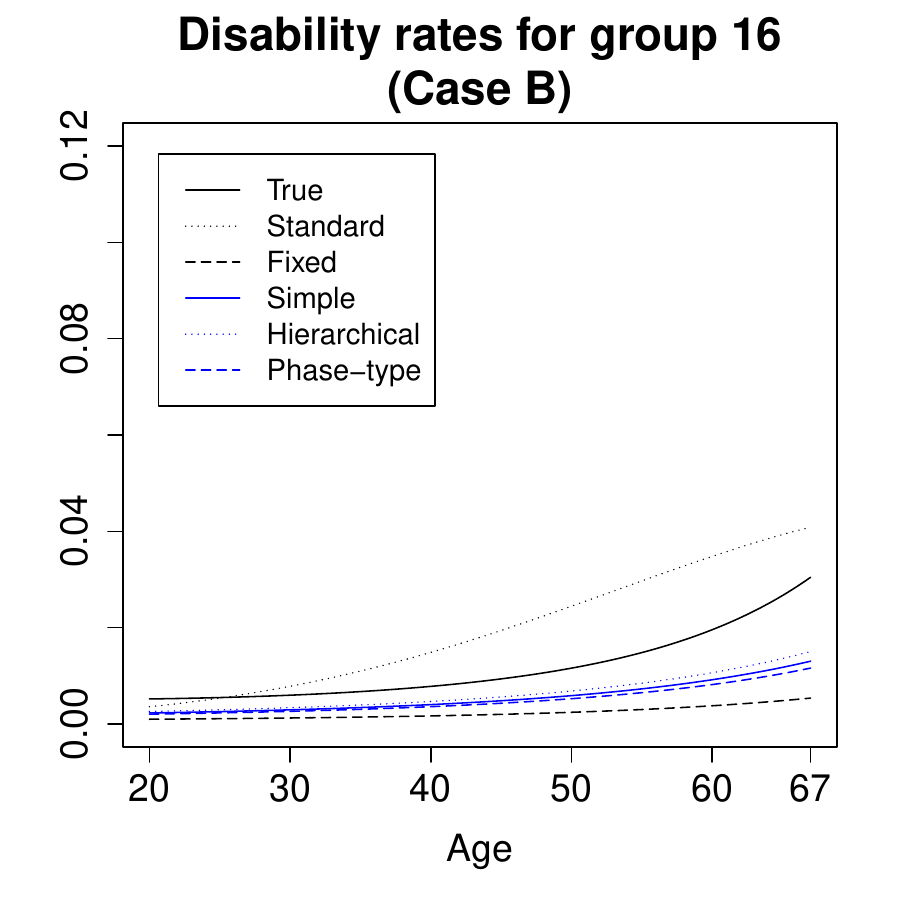}
\includegraphics[width=0.44\textwidth]{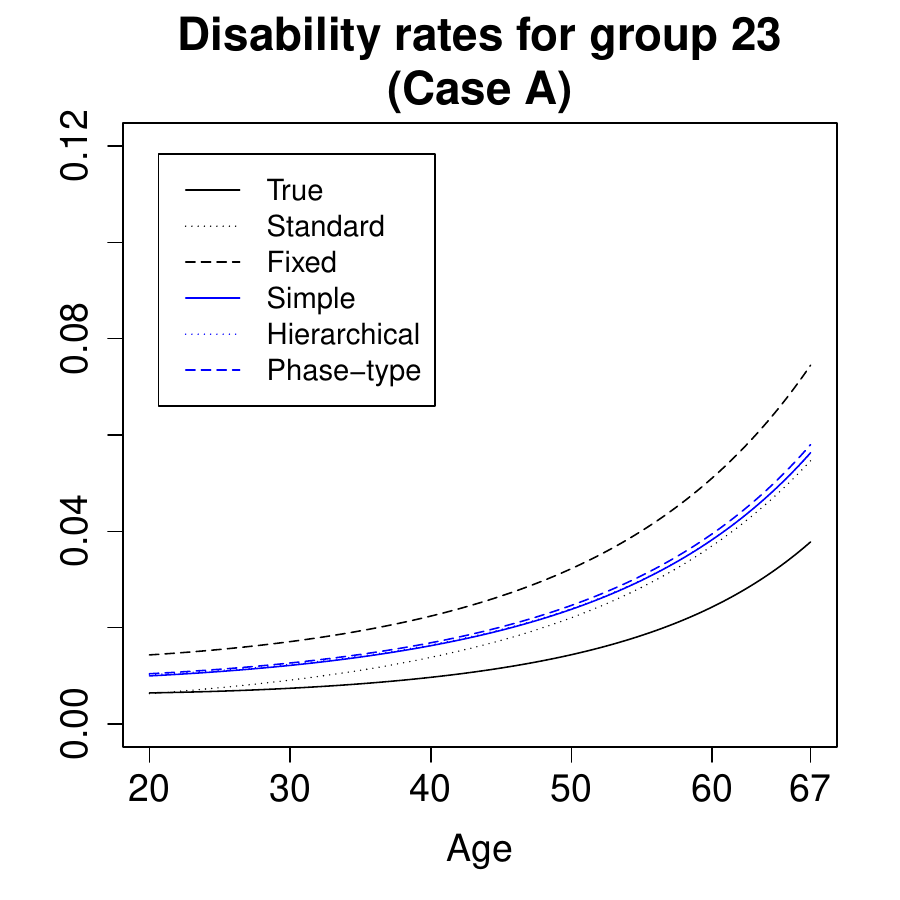}
\includegraphics[width=0.44\textwidth]{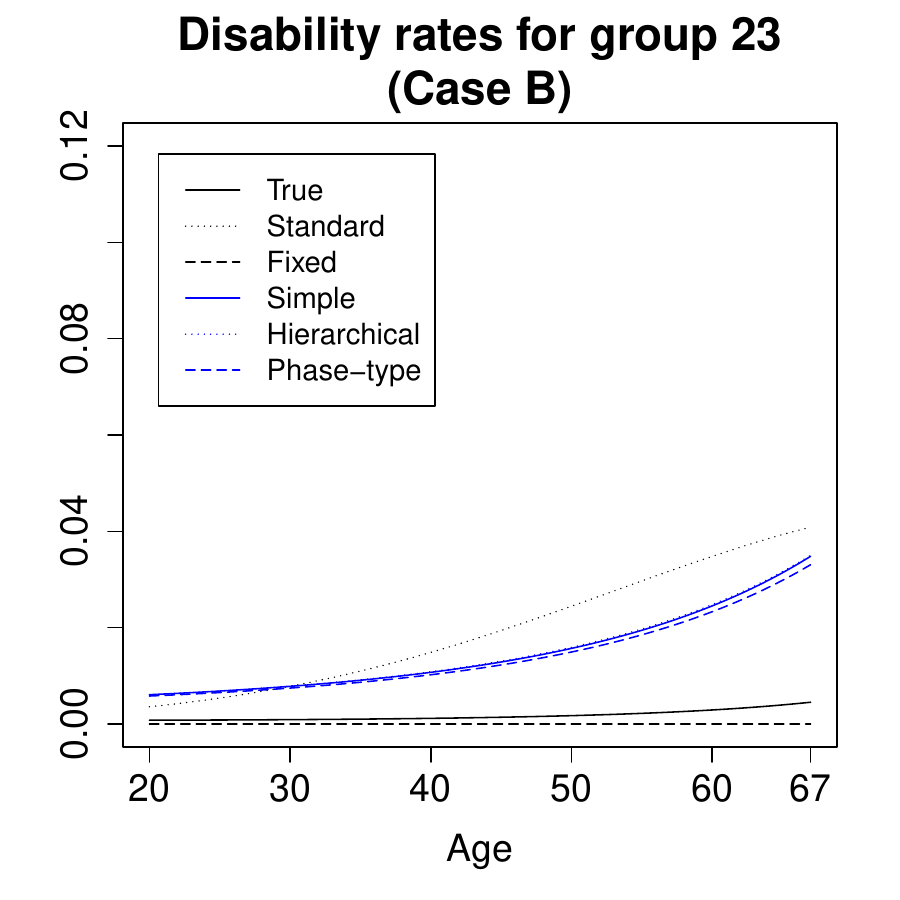}
\includegraphics[width=0.44\textwidth]{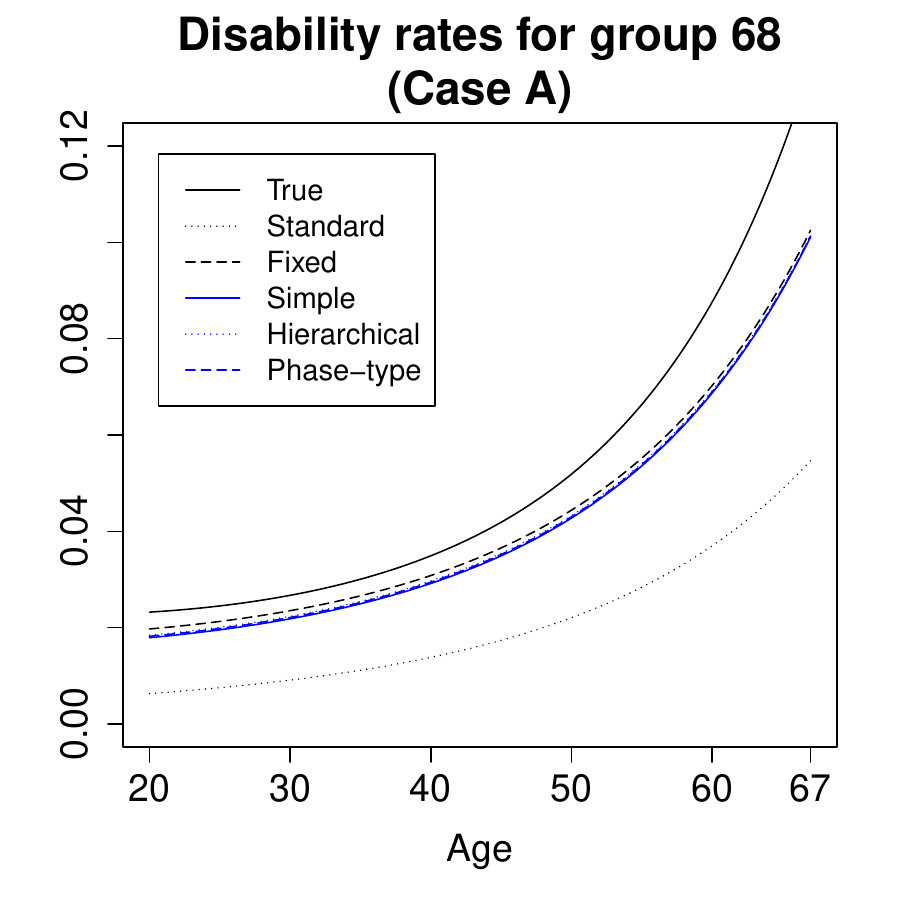}
\includegraphics[width=0.44\textwidth]{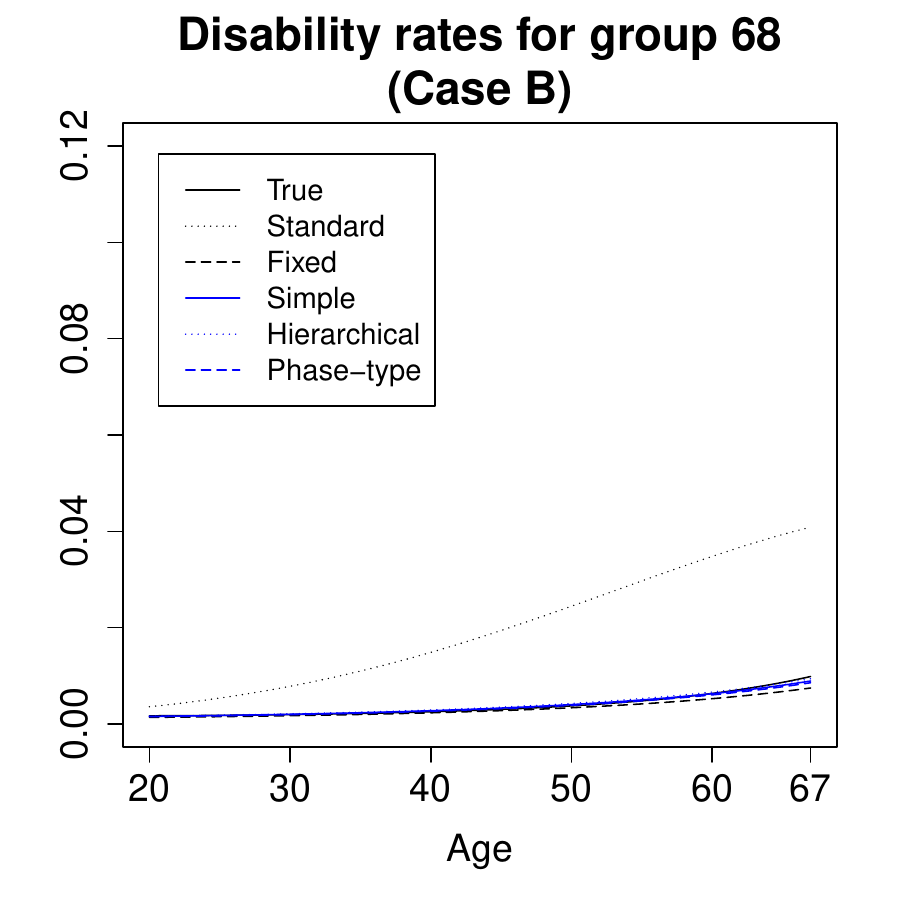}
\caption{Estimates of disability rates for the three selected groups in Case A and Case B.}
\label{fig:drates_s22_gamma_ind}
\end{figure}

\begin{figure}[!htbp]
\centering
\includegraphics[width=0.44\textwidth]{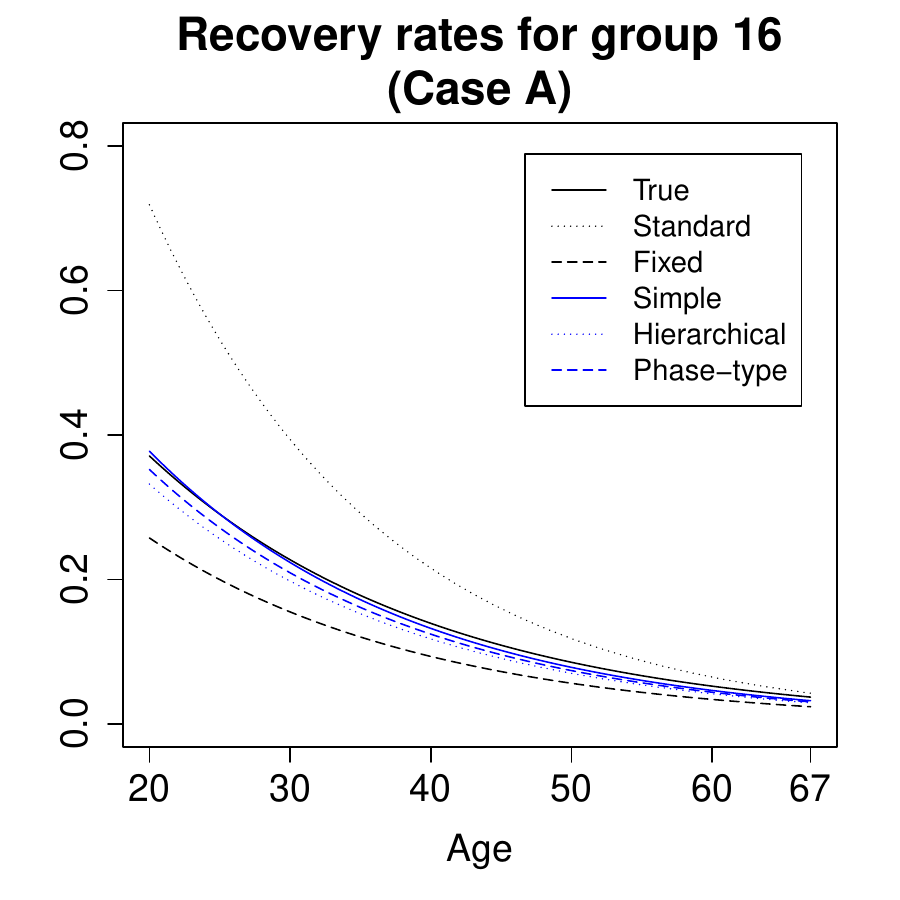}
\includegraphics[width=0.44\textwidth]{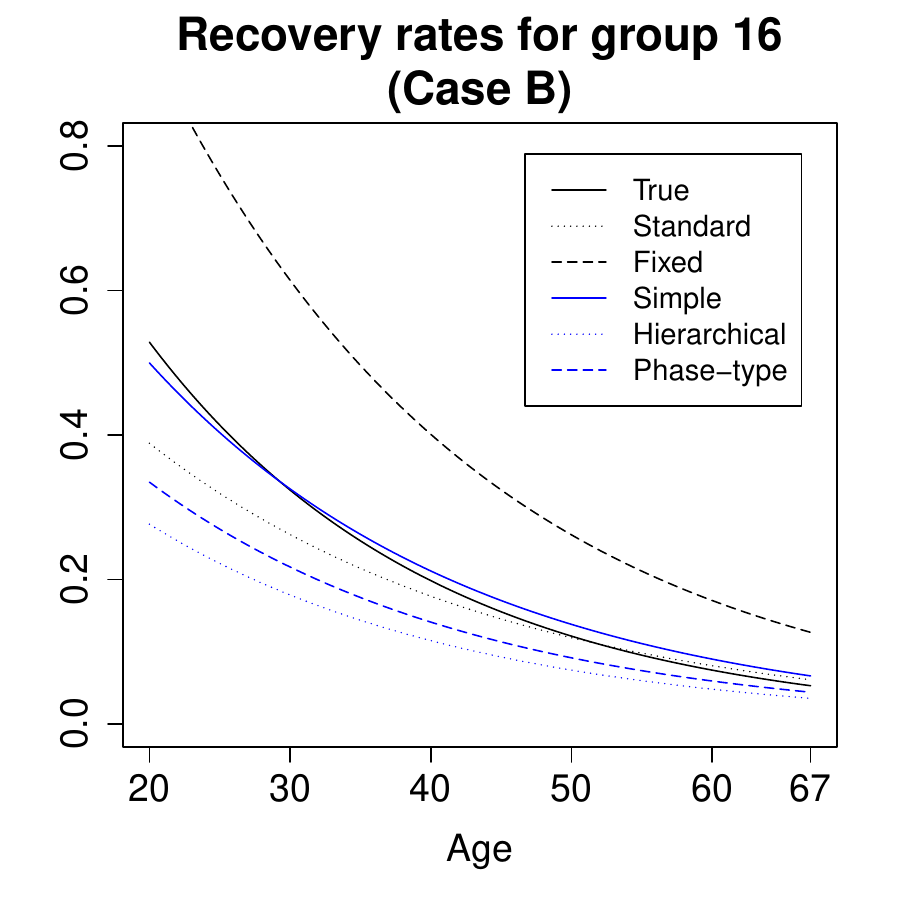}
\includegraphics[width=0.44\textwidth]{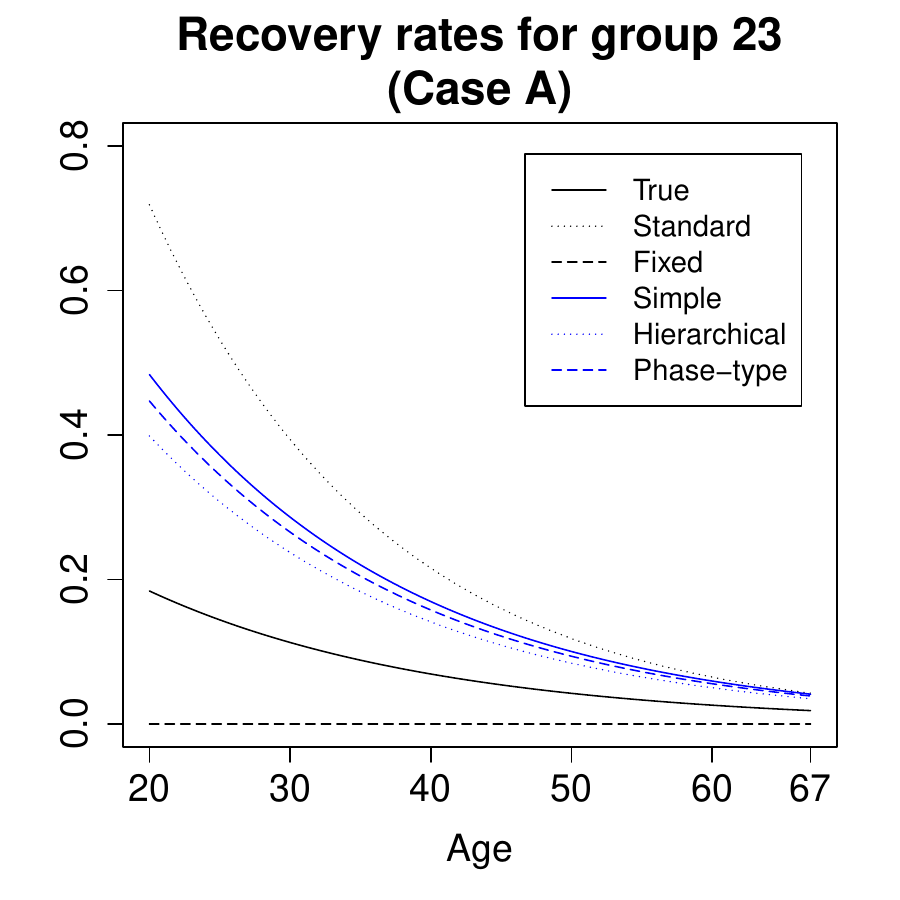}
\includegraphics[width=0.44\textwidth]{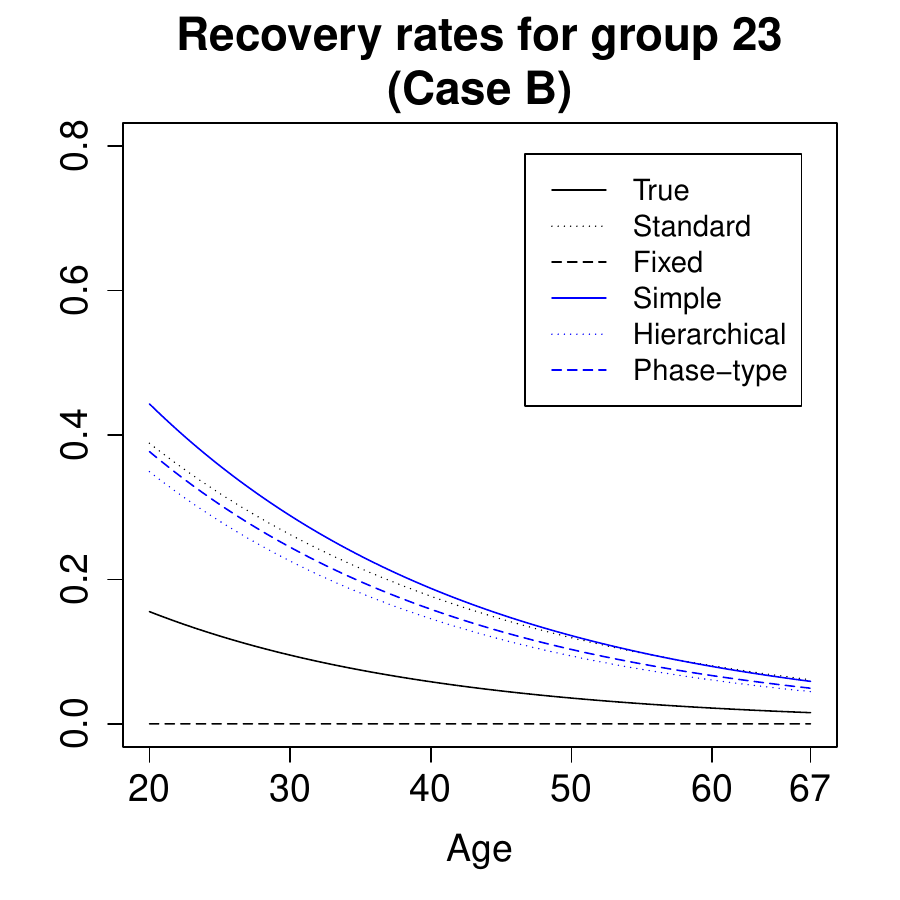}
\includegraphics[width=0.44\textwidth]{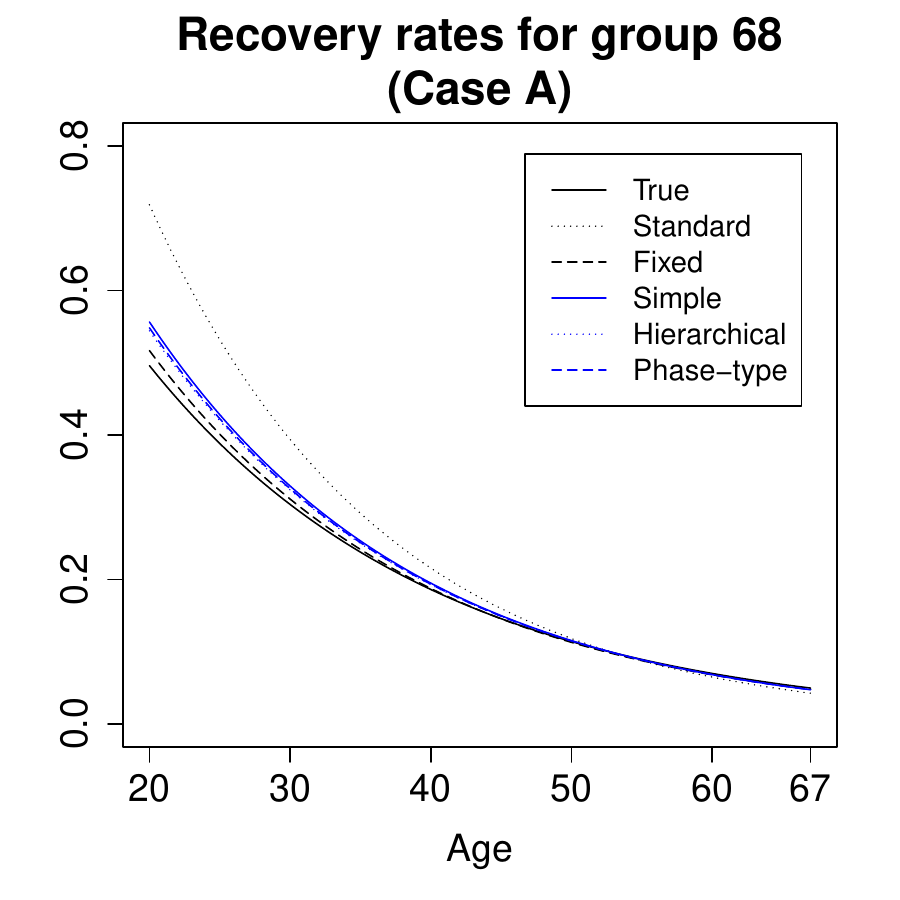}
\includegraphics[width=0.44\textwidth]{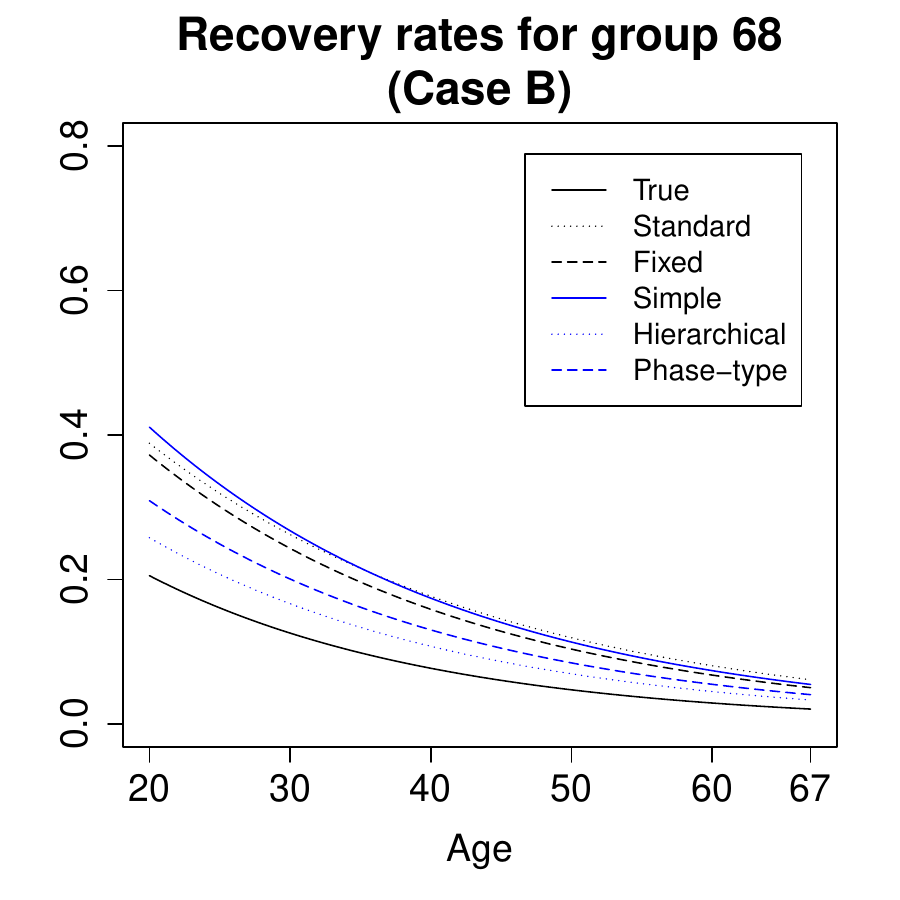}
\caption{Estimates of recovery rates for the three selected groups in Case A and Case B.}
\label{fig:rrates_s22_gamma_ind}
\end{figure}

For Case A, the estimates of the disability rates can be found in Figure~\ref{fig:drates_s22_gamma_ind}, while the corresponding recovery rates are illustrated in Figure~\ref{fig:rrates_s22_gamma_ind}. For the disability rates (Figure~\ref{fig:drates_s22_gamma_ind}), while the estimates of group $16$ are similar across all models considered, for the other two groups ($23$ and $68$), the three shrinkage fits outperform with similar results the standard and fixed effect fits. The superiority of the shrinkage estimates becomes more evident in the recovery rates (Figure~\ref{fig:rrates_s22_gamma_ind}), where they result in closer approximations to the true rates. This is particularly clear for group $23$, where the fixed effect model estimates a group effect of zero. We also observe the downside of the standard model, which yields identical estimates for all groups. 
{
As an additional visual aid to understand the results across all groups, Figure~\ref{fig:scatter_s22_gamma_ind} presents scatterplots of the (unscaled) estimated group effects under the different models, alongside the corresponding scatterplot of the (true) simulated ones for reference. For the simple model, we observe that several points align along a horizontal line. These correspond to groups without observed disabilities, resulting in absolutely no data on these group's particular recovery rate. In such a case, `full shrinkage' occurs: the estimated group effect essentially coincides with the mean of the random effect. This phenomenon appears not only for the simple model, but also for the hierarchical and phase-type model; for instance, for the former, the estimated value of $\nu$ is $40.525$, indicating a very small correlation. While the scatterplots are not directly comparable to the simulated values due to scaling, they provide a useful visual representation of the type of dependence each model is able to generate.
\begin{figure}[!htbp]
	\centering
		\includegraphics[width=0.44\textwidth]{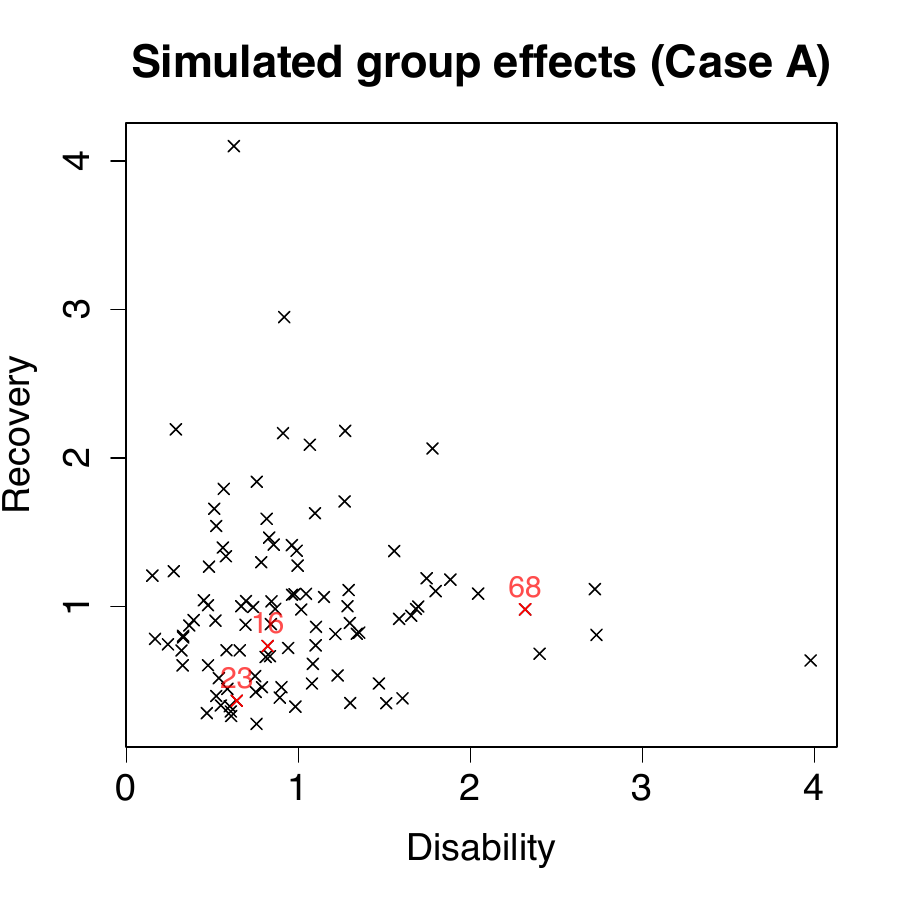}
		\includegraphics[width=0.44\textwidth]{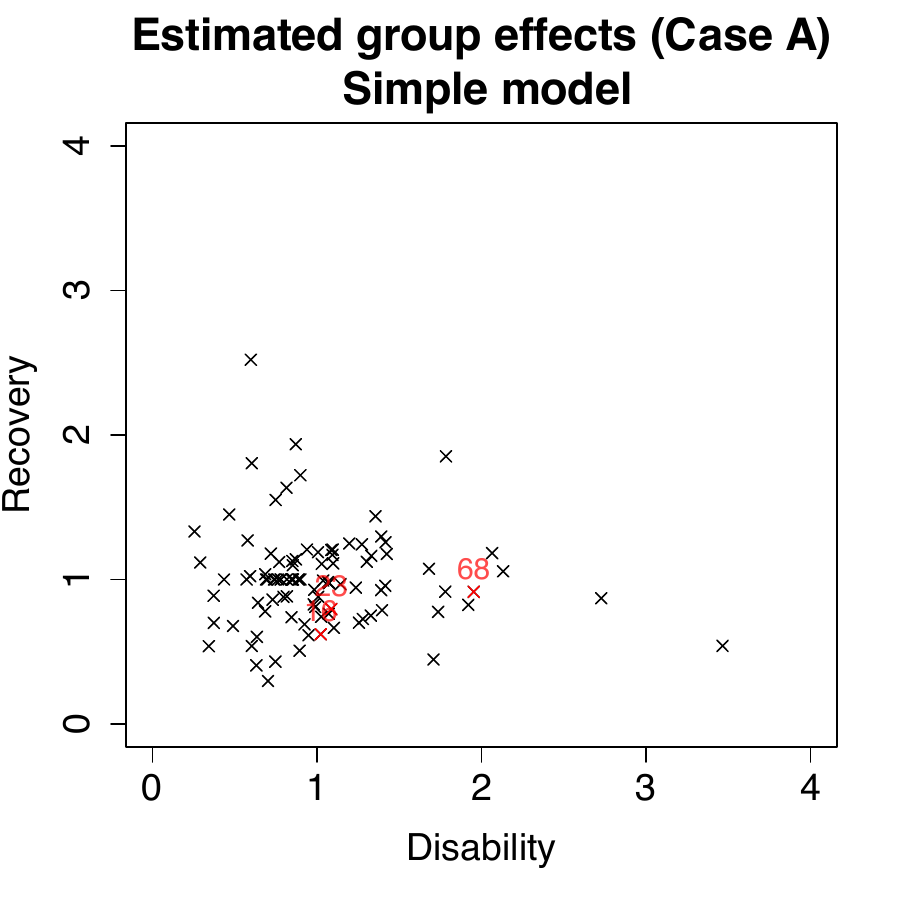}
		\includegraphics[width=0.44\textwidth]{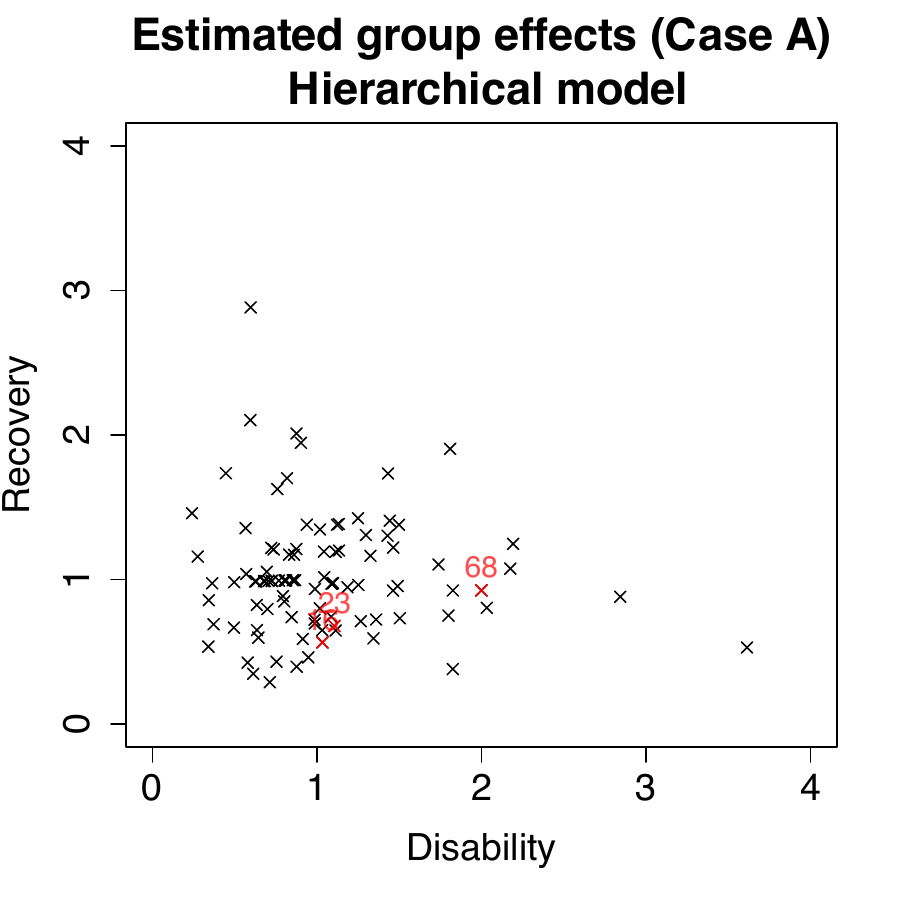}
		\includegraphics[width=0.44\textwidth]{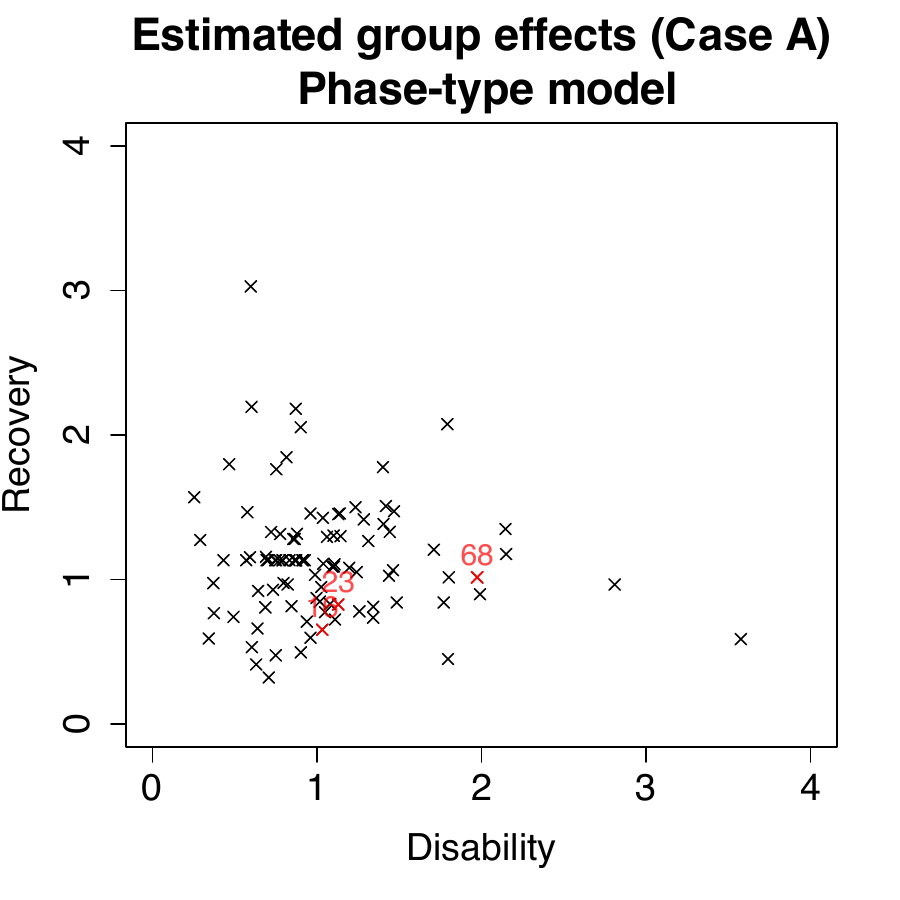}
\caption{Scatterplots of simulated and estimated group effects for Case A.}
\label{fig:scatter_s22_gamma_ind}
\end{figure}
}

Moving on to Case B, estimates of the disability rates are depicted in Figure~\ref{fig:drates_s22_gamma_ind}, and a similar conclusion can be drawn to that of Case A, with the shrinkage estimates showing a clear advantage over the standard and fixed effect fits. However, the primary distinction among the shrinkage estimates is observable in the recovery rates, presented in Figure~\ref{fig:rrates_s22_gamma_ind}. More specifically, for groups $16$ and $68$, the models capable of capturing positive correlation -- namely, the hierarchical and phase-type models -- provide closer approximations to the true disability rates. This leads us to conclude that by allowing dependency between latent group effects, the estimates of recovery and disability rates can mutually benefit.
{
This is confirmed by Figure~\ref{fig:scatter_s22_mgamma_pos}, which shows scatterplots of the simulated and estimated group effects. We observe that the hierarchical and phase-type models now produce estimates for groups without observed disabilities aligned along a positively sloped line, as they are able to learn the positive dependence structure. This is in stark contrast to the simple model, which remains unable to do so, resulting in points along a horizontal line.

\begin{figure}[!htbp]
	\centering
		\includegraphics[width=0.44\textwidth]{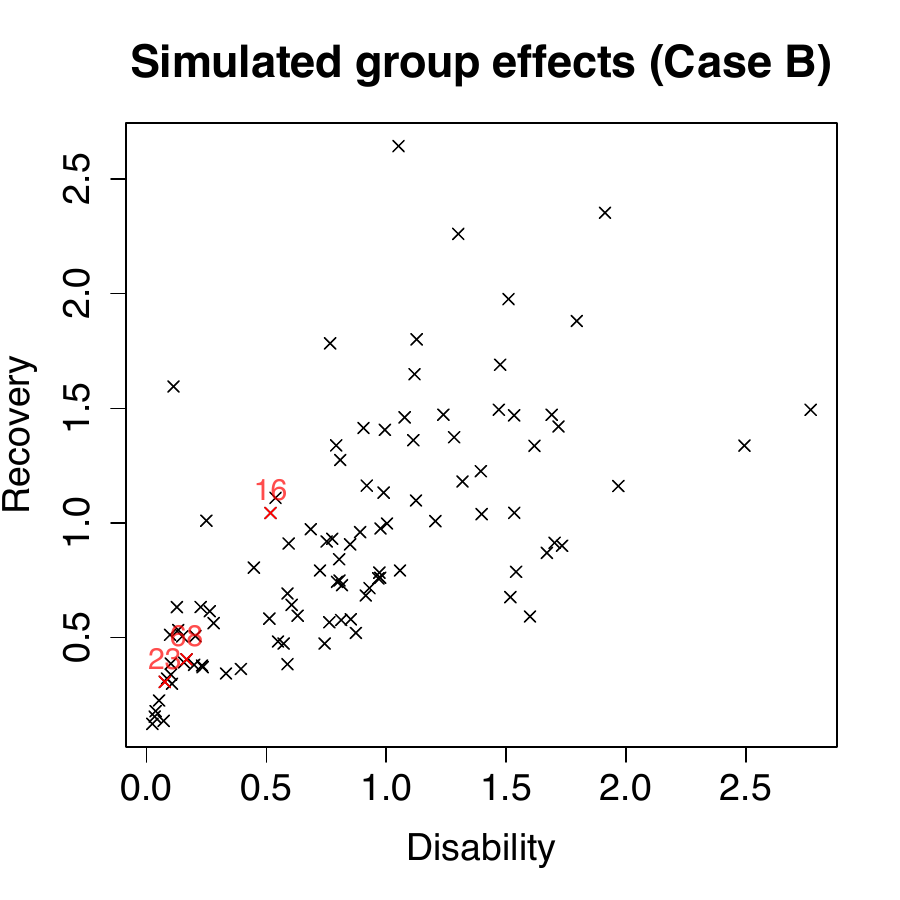}
		\includegraphics[width=0.44\textwidth]{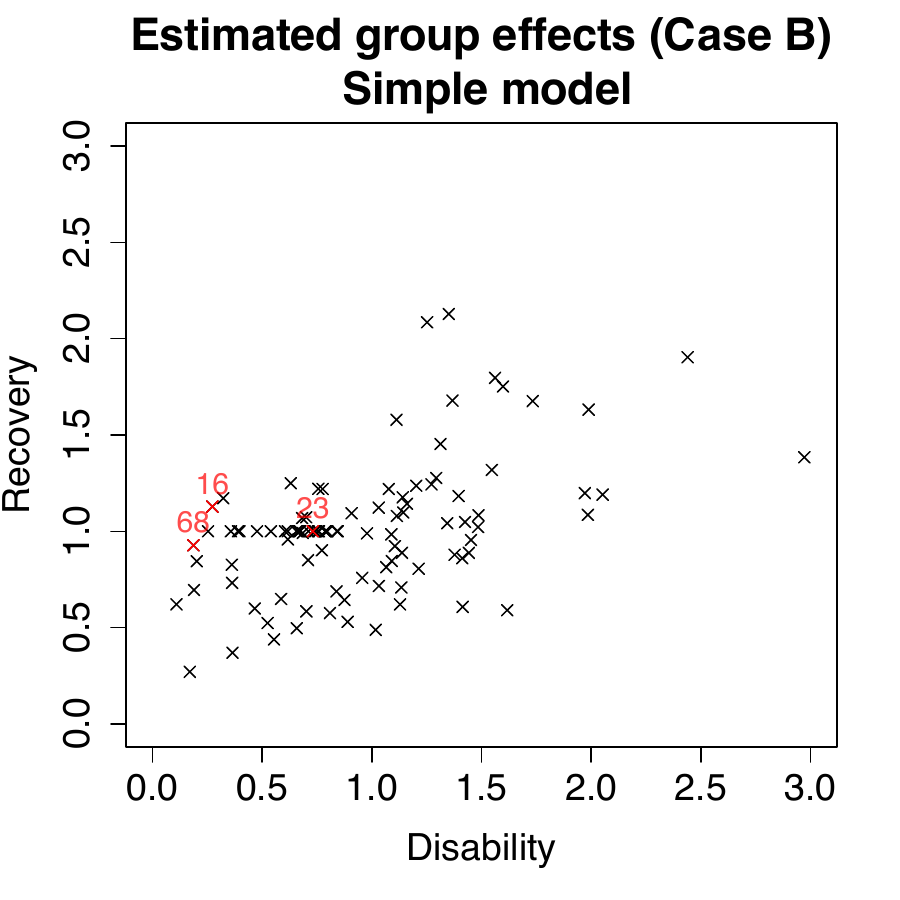}
		\includegraphics[width=0.44\textwidth]{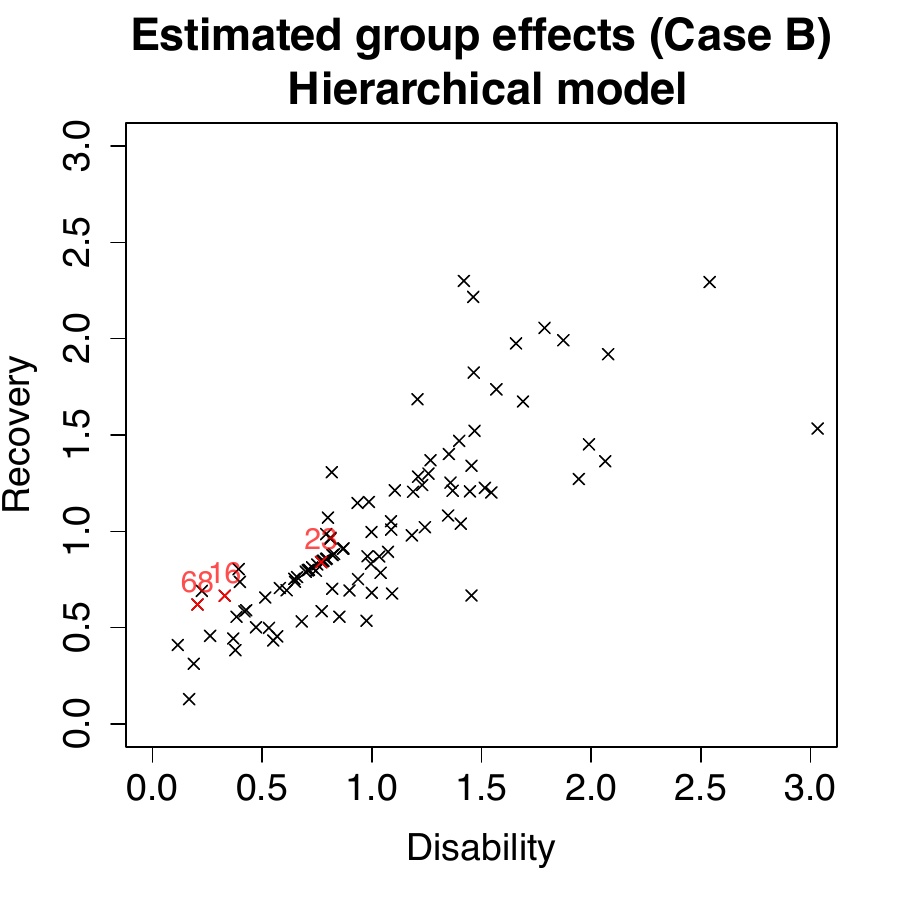}
		\includegraphics[width=0.44\textwidth]{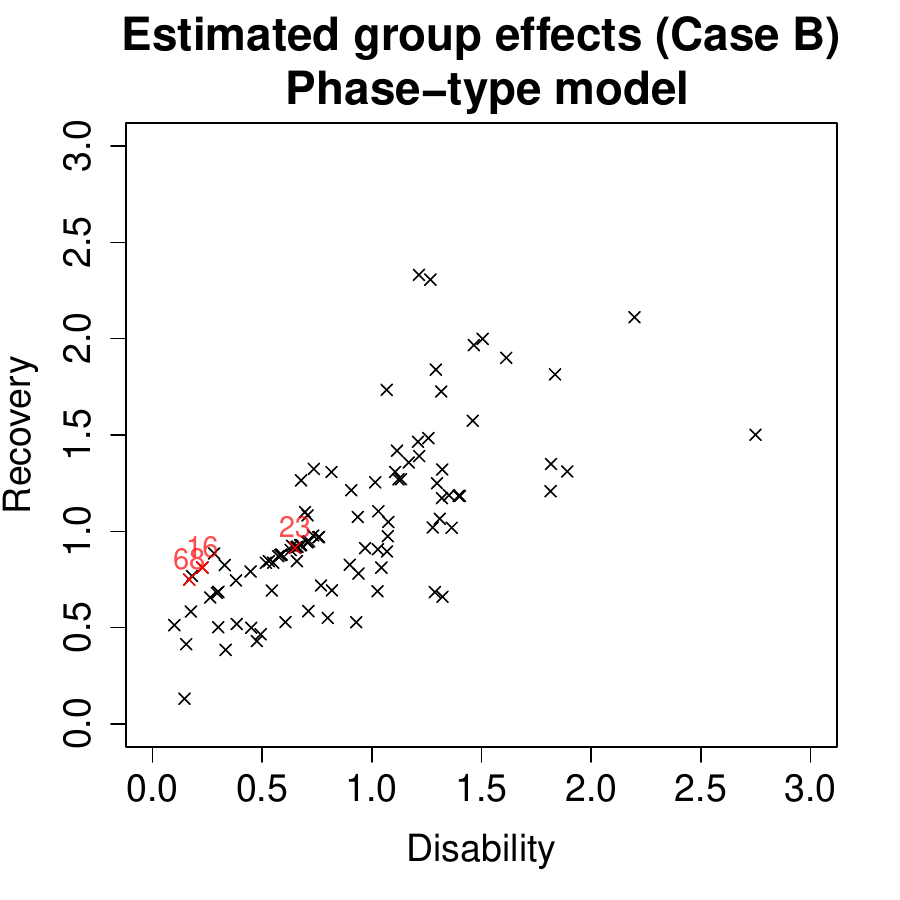}
\caption{Scatterplots of simulated and estimated group effects for Case B.}
\label{fig:scatter_s22_mgamma_pos}
\end{figure} 

Lastly, for Case C, we limit ourselves to presenting the scatterplots of the simulated and estimated group effects, Figure~\ref{fig:scatter_s22_mgamma_neg}, where we observe that the phase-type model is able to recover the negative dependence of the group effects. This is something the simple model cannot achieve. Jointly with the results from Cases A and B, this showcases the flexibility of the phase-type model to capture a wide range of dependence structures.

\begin{figure}[!htbp]
	\centering
		\includegraphics[width=0.44\textwidth]{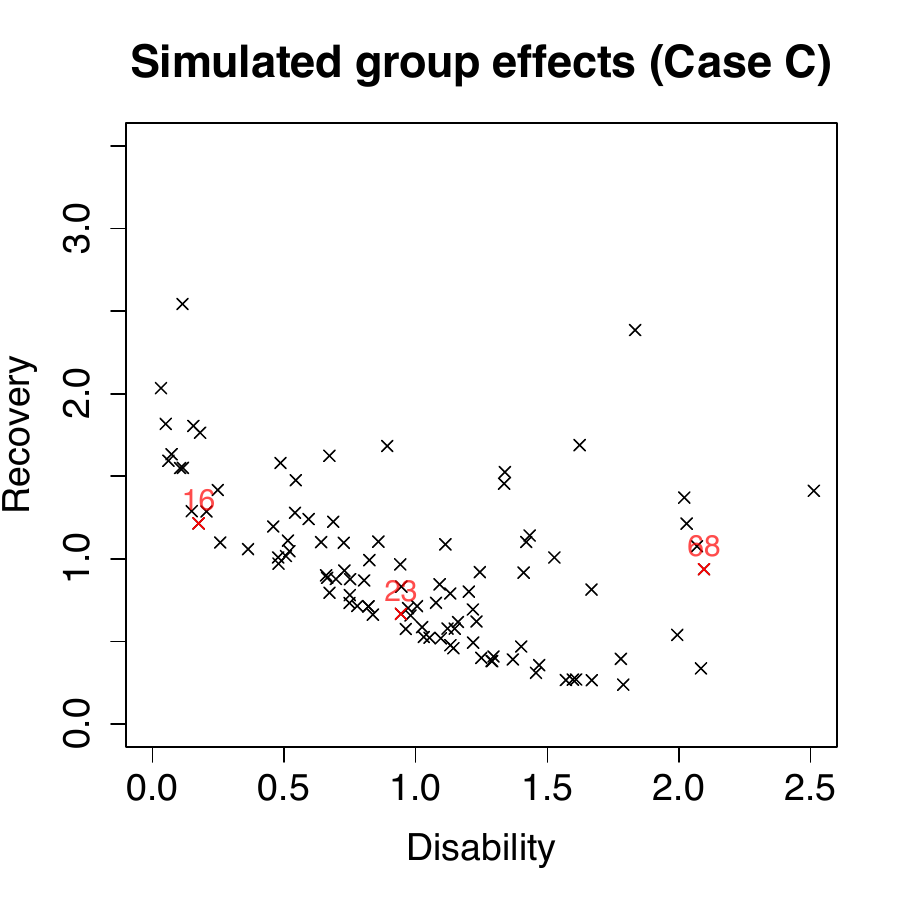}
		\includegraphics[width=0.44\textwidth]{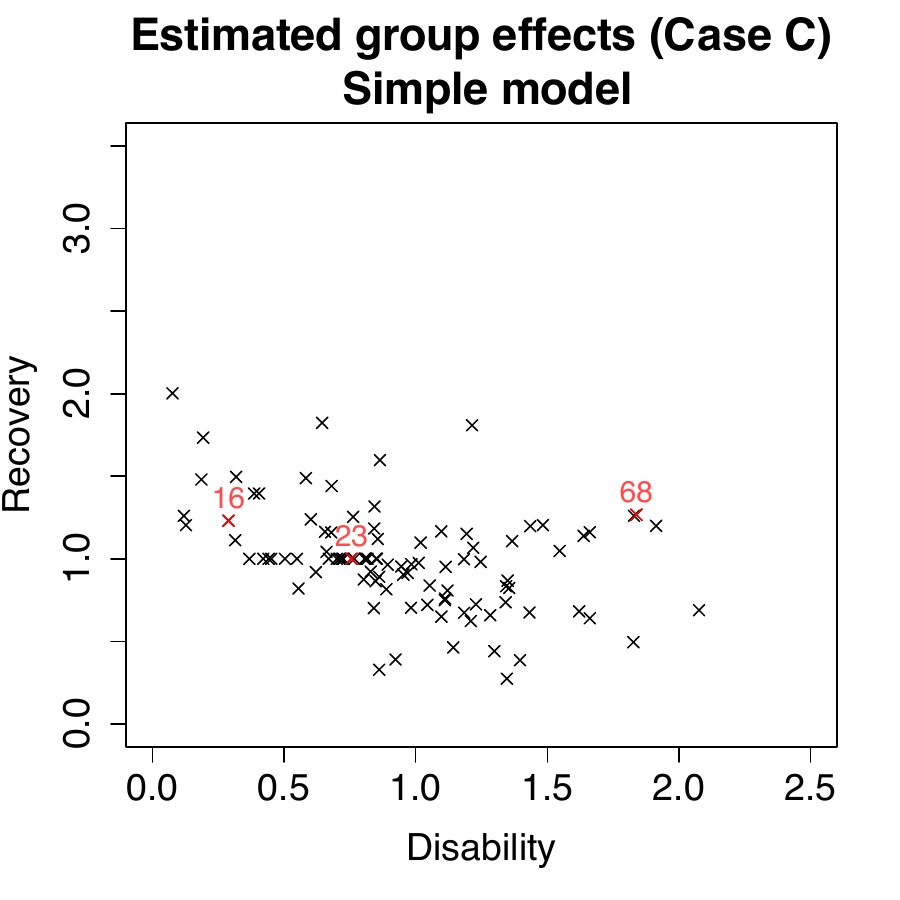}
		\includegraphics[width=0.44\textwidth]{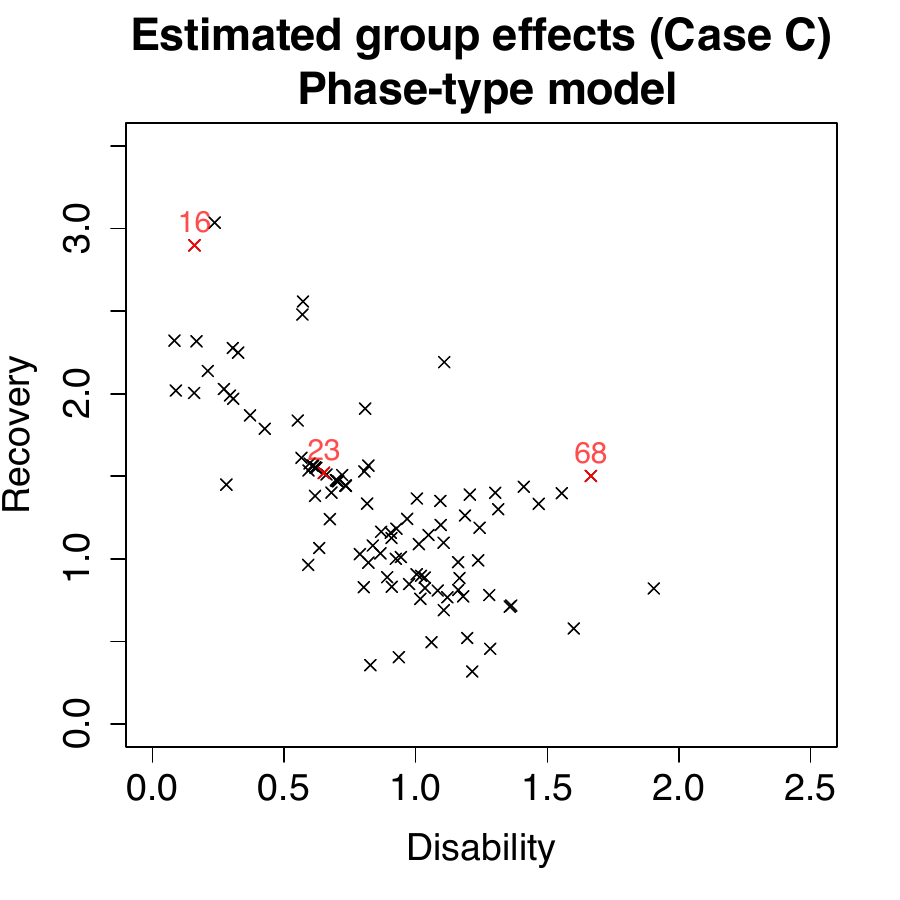}
\caption{Scatterplots of simulated and estimated group effects for Case C.}
\label{fig:scatter_s22_mgamma_neg}
\end{figure} 

}

We conclude this part of the analysis by presenting the loglikelihoods of the fitted mixed Poisson regression models. These loglikelihoods serve as statistical means to further assess the quality of the fits and can be found in Table~\ref{tab:loglik}. The observations from the table are consistent with what one could expect. In Case A, the simple model outperforms others, with the phase-type model following closely. This result aligns with the nature of the models and the data considered for this case. For Case B, the hierarchical model yields the highest loglikelihood, followed closely by the phase-type model. This outcome is somewhat anticipated since these two models are capable of capturing the positive correlation present in the mixing vector, aligning well with the underlying structure. Finally, for Case C, the phase-type model takes the lead, surpassing the { simple model}. This is hardly surprising, as it is the only model that can capture a negative correlation. Overall, the phase-type model proves to be very competitive across all cases. Its flexibility to capture different dependence structures highlight its potential as a robust tool for experience rating modeling.

\begin{table}[!htbp]
\centering
\begin{tabular}{l |c c c }
\hline
	Model & Case A   & Case B   & Case C  \\
 \hline 
Simple & $-3924.66$   & $-3816.12 $   & $-3639.01$   \\
Hierarchical & {$-3926.73$}   & {$-3803.57$}  & {---} \\
Phase-type        & $-3925.29$   & $-3805.56$  & $-3627.30$ \\
\hline
\end{tabular}
\caption{Loglikelihoods of the fitted mixed Poisson regression models. }
\label{tab:loglik}
\end{table}

\subsection{Predictive performance} \label{subsec:pred}

To undertake a comprehensive analysis of the predictive performance of the models considered, we carried out $100$ iterations of the simulation of individual disabilities and recoveries, keeping the group sizes, initial ages, and parameters, including the simulated group effects (for each case), fixed. Specifically, this results in $300$ datasets across the three cases (A, B, and C). For each dataset, we fitted the three models under consideration, allowing us to capture the nuances and variabilities in performance that may arise from different portfolio realizations.

\begin{figure}[!htbp]
\centering
\includegraphics[width=0.9\textwidth]{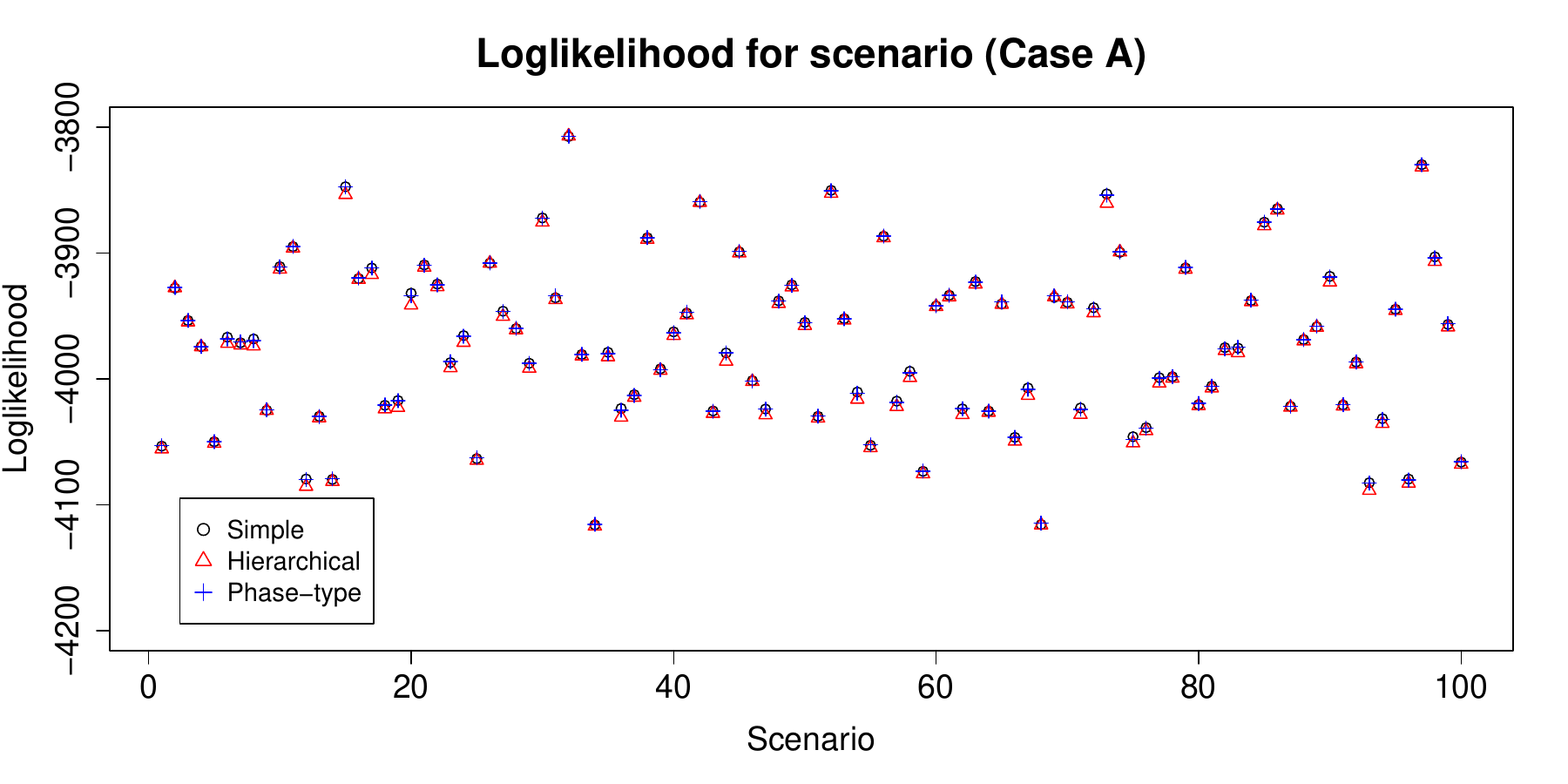}
\includegraphics[width=0.9\textwidth]{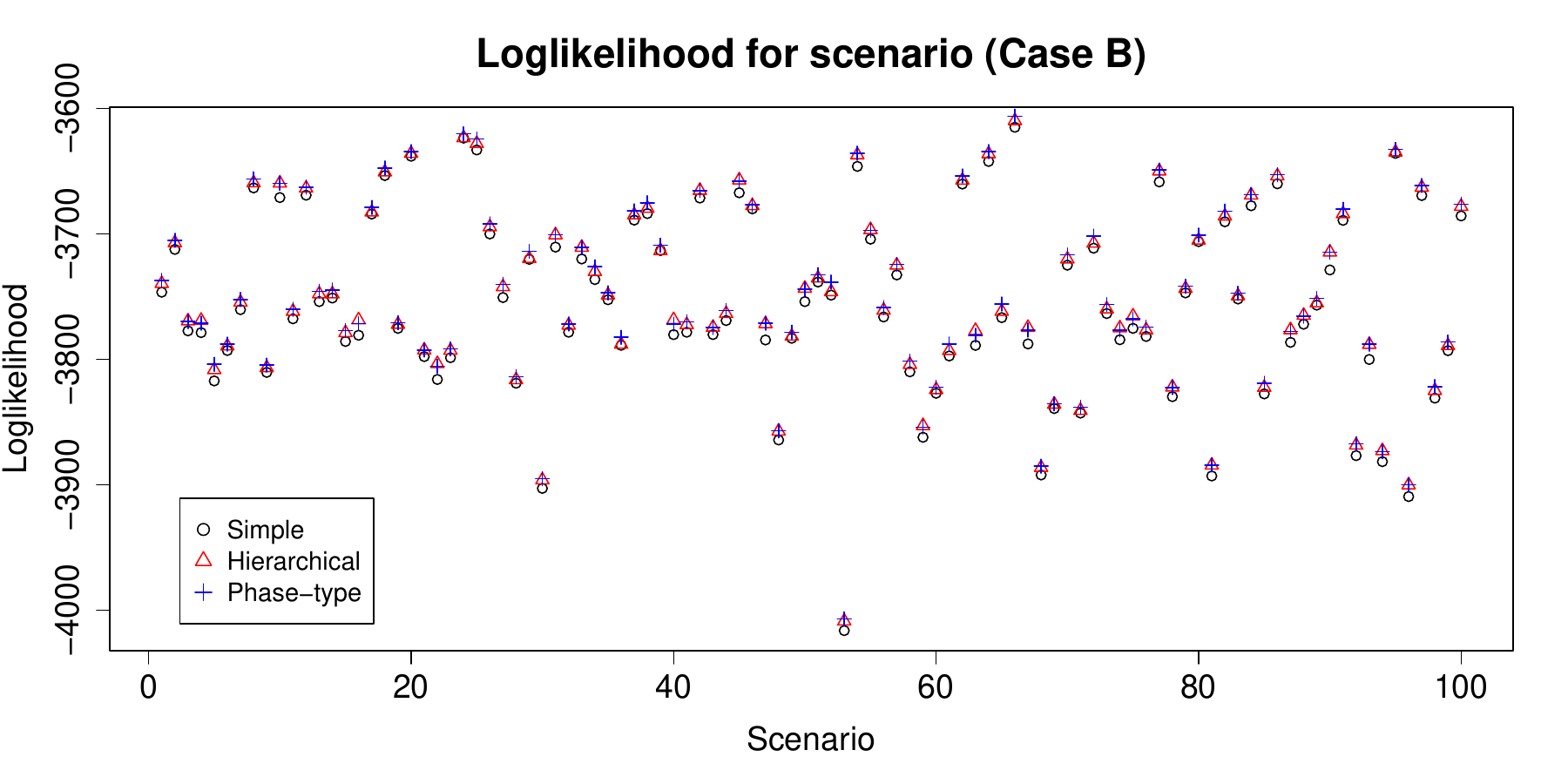}
\includegraphics[width=0.9\textwidth]{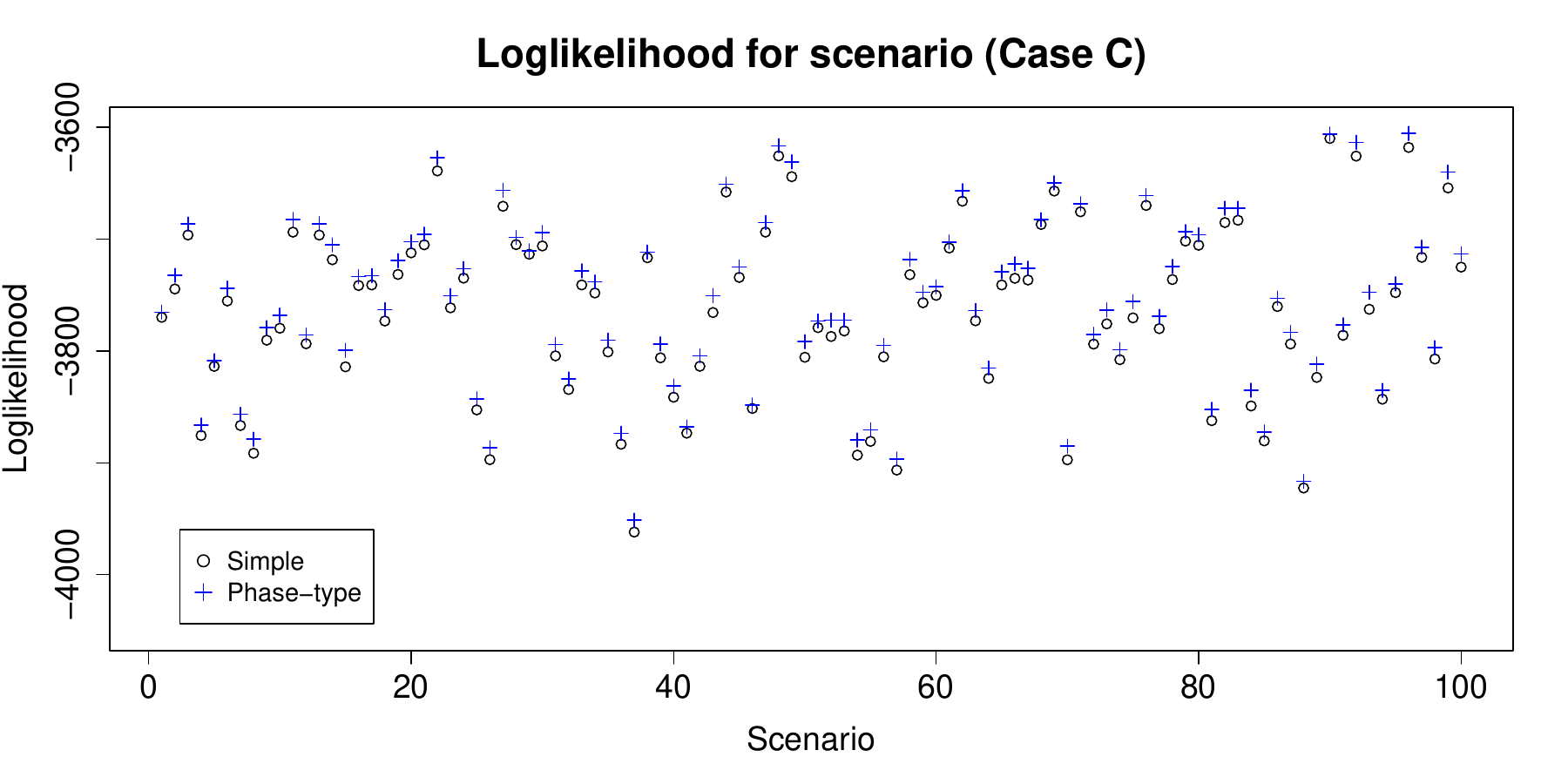}
\caption{Loglikelihoods of the fitted mixed Poisson regression models for $100$ iterations of the simulation study.}
\label{fig:ll}
\end{figure}

Figure~\ref{fig:ll} illustrates the loglikelihoods obtained in each of the three cases of interest across all iterations of the simulation study. The visual representation enables a more nuanced understanding of the models' behavior. More specifically, and consistent with our earlier observations, in Case A, the simple and phase-type models excel in all iterations. In contrast, the hierarchical model falls behind, which is expected given that this model can only describe positive correlations in the random effects. Nevertheless, Case B showcases the efficiency of the hierarchical model in capturing positive correlation, with the phase-type model also demonstrating strong performance and even outperforming the hierarchical model in several iterations.  Finally, in Case C, the phase-type model's ability to handle negative correlation sets it apart from the { simple model}, leading to the highest loglikelihoods for all iterations.

Overall, this extensive simulation study reaffirms that the phase-type model is the most versatile among the considered models. 

To evaluate the predictive performance of the models, we carry out premium calculations for a disability rate of one monetary unit until age $67$ under the equivalence principle using both the true parameters and the various estimates. A fourth order Runge-Kutta method with a step length of one month is employed. Premiums are first calculated for integer initial ages, and then linear interpolation is used to determine the premium for non-integer initial ages. Throughout, we use $r \equiv 0.01$, corresponding to a continuously compounded interest rate of $1 \%$. The comparative analysis is summarized in Table~\ref{tab:rmse}, which presents both the mean absolute errors (MAE's) and the root mean square errors (RMSE's) for each model under consideration and aggregated across groups and insured. 


\begin{table}[!htbp]
\centering
\begin{tabular}{l |r r r | r r r}
\hline
\multirow{2}{*}{Model} & \multicolumn{3}{c|}{RMSE} & \multicolumn{3}{c}{MAE} \\
\cline{2-7}
	& Case A   & Case B   & Case C & Case A   & Case B   & Case C \\
 \hline 
 Standard & $15.913$   & $14.989$   & $21.107$  & $10.450$ & $10.211$ & $15.158$ \\
 Fixed & $7.435$   & $6.511 $   & $7.887$  & $3.538$ & $3.457$ & $3.876$  \\
Simple & $5.538$   & $4.614 $   & $6.137$  & $3.131$ & $3.027$ & $3.634$  \\
Hierarchical & {$5.528$}   & {$4.456 $}   & {---} & {$3.115$} & {$2.940 $} & {---}  \\ 
Phase-type & $5.512$   & $4.509$  & $6.006$ & $3.114$ & $2.974$ & $3.517$  \\
\hline
\end{tabular}
\caption{RMSE's and MAE's for the total equivalence premium (in hundredths).}
\label{tab:rmse}
\end{table}

The findings mirror our previous insights, reinforcing the supremacy of the shrinkage models over the standard and fixed effect approaches. 
{
The phase-type model exhibits the most favorable performance in Cases A and C, while in Case B the hierarchical model takes the lead, with the phase-type model following closely. Overall, this underlines the broader flexibility of the phase-type specification.

Finally, it is essential to note that the hierarchical model's assumption of equal marginal variances is quite limiting. For example, when considering a scenario where $\Theta_{ai}$ follows a $\text{Gamma}(10/6, 10/6)$-distribution and $\Theta_{ia}$ follows a $\text{Gamma}(10/3, 10/3)$-distribution, coupled by the same Clayton copula with Kendall's $\rho_\tau = 0.5$, even the simple model can outperform the hierarchical model, with an improvement (reduction) in RMSE of about $0.03$. This showcases that, beyond capturing dependence, accounting for differences in marginal variances is also crucial for predictive performance.
}

\section{Conclusion}\label{sec:conclusion}

This paper introduces two novel and tractable specifications for mixed Poisson regression: the hierarchical and the phase-type. For each model, an estimation method via an expectation-maximization algorithm has been derived. Further, following~\cite{furrer2019experience}, we have showcased how these methods might be applied for shrinkage estimation of group effects in the context of disability insurance. The simulation study confirms the effectiveness of the phase-type model for such tasks. Our results reval that, unlike the simple or hierarchical models, which are constrained by their inherent dependence structure assumptions, the phase-type model is a robust and verstaile tool across various realistic cases. We should like to stress that the potential of both the hierarchical as well as the phase-type specification extends beyond the immediate context, and both models could find application in other areas. For instance, it would be interesting to study their performance for frequency modeling in the context of non-life insurance.

\end{document}